\documentclass[aap,MSNbibl,nameyear,dvips]{arximspdf}
\usepackage{graphicx}

\doi{10.1214/11-AAP767}
\volume{22}
\issue{1}
\pubyear{2012}
\firstpage{239}
\lastpage{284}

\makeatletter

\def\bsuffix #1{#1}

\newcommand{\trace}{\operatorname{tr}}
\newcommand{\espp}[1]{E_P[#1]}
\newcommand{\espalt}[2]{E_{#1}[#2]}
\newcommand{\mfnct}{\eta}
\newcommand{\hatp}{\pi}
\newcommand{\lrpair}{(\hatp,\mfnct)}
\newcommand{\reals}{\mathbb{R}}
\newcommand{\borel}{\mathcal{B}}
\newcommand{\eps}{\varepsilon}
\newcommand{\ups}{\Upsilon}
\newcommand{\conv}{\rightarrow}
\newcommand{\der}[1]{\dot{#1}}
\newcommand{\dder}[1]{\ddot{#1}}
\newcommand{\pit}{{\pi^T}}
\newcommand{\etat}{{\eta^T}}
\newtheorem{lemma}{Lemma}
\newproclaim{defin}[lemma]{Definition}
\newtheorem{proposition}[lemma]{Proposition}
\newtheorem{corollary}[lemma]{Corollary}
\newtheorem{theorem}[lemma]{Theorem}
\newproclaim{example}[lemma]{Example}
\newproclaim{remark}[lemma]{Remark}
\newproclaim{assum}[lemma]{Assumption}

\newcommand{\cale}{\mathcal{E}}
\newcommand{\calf}{\mathcal{F}}

\newcommand{\eqref}[1]{(\ref{#1})}
\makeatother

\begin{document}
\begin{frontmatter}

\title{Portfolios and risk premia for the long run\thanksref{T1}}
\runtitle{Portfolios and risk premia for the long run}

\thankstext{T1}{Supported in part by NSF (DMS-053239 and DMS-0807994), SFI
(07/MI/008, 07/SK/M1189, 08/SRC/FMC1389) and the European Commission
(RG-248896).}

\begin{aug}
\author[A]{\fnms{Paolo} \snm{Guasoni}\corref{}\ead[label=e1]{guasoni@bu.edu}}
\and
\author[B]{\fnms{Scott} \snm{Robertson}\ead[label=e2]{scottrob@andrew.cmu.edu}}
\runauthor{P. Guasoni and S. Robertson}
\affiliation{Boston University, Dublin City University and Carnegie
Mellon University}
\address[A]{Department of Mathematics and Statistics\\
Boston University\\
111 Cummington st.\\
Boston, Massachusetts 02215\\
USA\\
and\\
School of Mathematical Sciences\\
Dublin City University\\
Glasnevin, Dublin 9\\
Ireland\\
\printead{e1}}
\address[B]{Department of Mathematical Sciences\\
Carnegie Mellon University\\
Wean Hall 6113\\
Pittsburgh, Pennsylvania 15213\\
USA\\
\printead{e2}}
\end{aug}

\received{\smonth{7} \syear{2009}}
\revised{\smonth{2} \syear{2011}}

%
\begin{abstract}
This paper develops a method to derive optimal portfolios and risk
premia explicitly in a general diffusion model for an investor with
power utility and a long horizon. The market has several risky
assets and is \mbox{potentially} incomplete. Investment opportunities are
driven by, and partially correlated with, state variables which
follow an autonomous diffusion. The framework nests models of
stochastic interest rates, return predictability, stochastic
volatility and correlation risk.

In models with several assets and a single state variable, long-run
portfolios and risk premia admit explicit formulas up the solution
of an ordinary differential equation which characterizes the
principal eigenvalue of an elliptic operator. Multiple state
variables lead to a~quasilinear partial differential equation which
is solvable for many models of interest.

The paper derives the long-run optimal portfolio and the long-run
optimal pricing measures depending on relative risk aversion, as
well as their finite-horizon performance.
\end{abstract}

%


\end{frontmatter}

\section*{Introduction}\label{intr}
Long-run asymptotics are a powerful tool to obtain explicit formulas
in portfolio choice and derivatives pricing but their use has been
mostly ad hoc in the absence of general results. This paper
develops a~method to derive optimal portfolios and risk premia
explicitly in a general diffusion model for an investor with power
utility and in the limit of a long horizon. The market has several
risky assets and is potentially incomplete. Investment opportunities
are driven by, and partially correlated with, state variables that
follow an autonomous diffusion.

Investment and pricing problems share a reputation for mathematical
complexity. This common trait is not an accident; the central
message of duality theory\setcounter{footnote}{1}\footnote{See, for example,
\citet{MR844013}, \citet{MR912456}, \citet{MR1024460}, \citet
{MR1122311}, \citet{MR1722287}, as well as
\citet{MR2178034} for a setting similar to the one in this paper.} is
that the
two problems are indeed equivalent, as state-price densities are
proportional to the marginal utilities of optimal payoffs. In spite
of this conceptual equivalence, portfolio choice and derivatives
pricing have followed largely different strands of literature, each
of them with its own terminology.

The portfolio choice literature focuses on finding the
\textit{intertemporal hedging} component of optimal
portfolios.\footnote{\citet{kim1996dnp}, \citet{brennan2002daa},
\citet{wachter2002pac}, Munk and S{\o}ren\-sen (\citeyear{munk2004oca}), \citet{liu2007pss},
compute optimal portfolios explicitly for certain models.} Long-run
asymptotics have appeared in this literature under different names:
the \textit{risk-sensitive control}
approach,\footnote{
\citet{MR1358100}, \citeauthor{MR1372906} (\citeyear{MR1372906,MR1972534}),
\citeauthor{MR1675114}
(\citeyear{MR1675114,MR1790132}), \citeauthor{MR1802598}
(\citeyear{MR1802598,MR1910647}), \citet{MR1932164} and
\citet{MR1890061}, \citet{hata2005solving}, \citeauthor{MR2076543}
(\citeyear{MR2076543,MR2206349}).} turnpike
results\footnote{\citet{lelan}, \citet{hakanss}, \citet{MR736053},
\citet{MR1172445}, \citet{MR1629559}, \citet{MR1805320}, \citet
{backdybrog}.} and
large deviations
criteria\footnote{\citet{MR1968944}, \citet{follmeraaa}.} are all
efforts to achieve tractability by means of the long-run limit.

The derivatives pricing literature strives to identify martingale
measures that are optimal in the sense of the \textit{minimax
martingale measure} of \citet{MR1122311} or the \textit{least
favorable completion} of \citet{MR1089152}.
Power utility leads to the $q$-optimal measure
(\citep{MR2092922}, \citep{MR2116796}) which embeds several other martingale
measures;\footnote{Logarithmic utility leads to the minimal measure
used by \citet{MR1108430}. Exponential utility leads to the
minimal-entropy measure Grandits and Rheinl{\"a}n\-der (\citeyear{MR1920099}), \citet{MR1743972},
\citet{MR2211715}. Mean--variance hedging leads to the~%
va\-riance-optimal measure introduced by
\citeauthor{MR1143398} (\citeyear{MR1143398,MR1387633}).} it reduces to the minimal measure for
$q=0$, to the minimal entropy measure for $q=1$ and to the
variance-optimal measure for $q=2$.

The advantages of long-run asymptotics are their tractability and
accuracy. Long-run portfolios and risk premia are much simpler than
their finite-horizon counterparts and allow explicit expressions
even in cases in which the latter do not. In general, long-run policies
are identified by the quasilinear partial differential equation
\eqref{EmainvPDE} which admits explicit solutions in several
models of interest. In the case of a single state variable,
this equation reduces to an ordinary
differential equation, which
is furthermore linear if the state variable has a constant
correlation with excess returns. The accuracy of the long-run
approach stems from the bounds \eqref{cebou} which estimate the
duality gap at any horizon and hence the potential departure from
the unknown finite-horizon optimum. \textit{Long-run optimality} holds
(Definition \ref{lrdef}) when long-run policies are
approximately\vadjust{\goodbreak}
optimal over long horizons. The main result of this paper gives a
sufficient condition for long-run optimality in a general
multidimensional diffusion. Furthermore, this condition is sharp for
certain models and a calibration to the parameters estimated by
\citet{barberis2000ilr} shows that it is satisfied for reasonable
levels of risk aversion.

Two duality insights are central to our results. First,
the usual duality between payoffs and martingale densities extends
to their stochastic logarithms, which are portfolios and risk
premia.
Second, long-run asymptotics become easier in a duality context because
candidate long-run risk-premia yield an upper bound on the maximal
expected utility and vice versa. This observation allows to overcome
some difficulties arising in the verification theorems of the
\textit{risk-sensitive control} literature.

An important concept arising in long-run analysis is the
\textit{myopic probabili\-ty}; a~long-run investor with power utility
under the original probability beha\-ves like a logarithmic (or myopic)
investor under the myopic probability. This
probability plays an important role both for long-run analysis and for
finite-horizon bounds and its existence is crucial for the long-run
optimality result.\looseness=-1

The rest of the paper is organized as follows. Section \ref{secmodel}
describes the model in detail, introducing notation. Section \ref
{seclongrun} contains the
main result: a~general method to obtain long-run policies in closed
form. It also provides sufficient conditions, adapted from \citet
{MR2206349}, for the existence of solutions to the associated ergodic
Bellman equation. Section \ref{interpret} discusses the various
implications of these
results for portfolio choice and derivatives pricing, and the
connections with the
stochastic control and large deviations approaches. Section
\ref{secapp} derives long-run portfolios and risk premia in two
models of interest. The last one combines stochastic interest rates,
drifts and volatilities, and still admits simple closed form
solutions. Section \ref{secconc} concludes. All proofs are in the
appendices.\vspace*{-1pt}

\section{Model}\label{secmodel}
\subsection{Market}\label{sec1.1}
Consider a financial market with a risk-free asset $S^0$ and~$n$
risky assets $S=(S^1,\ldots,S^n)$. Investment opportunities (interest
rates, expected returns and covariances) depend on $k$
state variables $Y=(Y^1,\ldots,Y^k)$ which model their change over time,
\begin{eqnarray}
\label{rateeq}
\frac{dS^0_t}{S^0_t}&=&r(Y_t)\,dt,\\
\label{priceeq}
\frac{dS^i_t}{S^i_t}&=&r(Y_t) \,dt+dR^i_t, \qquad{1\le i\le n.}
\end{eqnarray}
Cumulative excess returns $R=(R^1,\ldots,R^n)$ and state variables
follow the diffusion
\begin{eqnarray}
\label{reteq}
dR^i_t&=&\mu_i(Y_t)\,dt+\sum_{j=1}^n \sigma_{ij}(Y_t) \,dZ^j_t,
\qquad{1\le i\le
n},\\
\label{stateeq}
{dY^i_t}&=& b_i(Y_t)\,dt+\sum_{j=1}^k a_{ij}(Y_t)\,dW^j_t, \qquad{1\le
i\le k},\\
\label{correq}
d\langle Z^i,W^j\rangle_t &=& \rho_{ij}(Y_t)\,dt, \qquad1\le i\le n,
1\le j\le k,
\end{eqnarray}
where $Z=(Z^1,\ldots,Z^n)$ and $W=(W^1,\ldots,W^k)$ are multivariate
Brownian motions. This setting provides a flexible framework that
nests most diffusion models in finance, including the models of
correlation risk considered by \citet{buraschi2010correlation} in
which $\rho$ is a function of a state variable.


The law of $(R,Y)$ determines the drifts $b, \mu$ and the
covariation matrices $\Sigma=\sigma\sigma'={d\langle
R,R\rangle_t}/{dt}$, $A=aa'={d\langle Y,Y\rangle_t}/{dt}$ and
$\ups=\sigma\rho a'={d\langle R,Y\rangle_t}/{dt}$, where the prime
sign denotes matrix transposition. By contrast, the matrices
$\sigma, a, \rho$ are identified only up to orthogonal transformations.
The market defined by \eqref{rateeq}--\eqref{correq} is in general
incomplete and
the covariance matrix
$\ups'\Sigma^{-1}\ups$ gauges the degree of incompleteness of the
market, highlighting two extremes:\vspace*{1pt} \textit{complete markets} for
$\ups'\Sigma^{-1}\ups= A$ and \textit{fully incomplete markets} for
$\ups=0$.

Let $E\subseteq\reals^k$ be an open connected
set. Denote by $C^{m}(E,\reals^{d})$
[resp., $C^{m,\gamma}(E,\allowbreak\reals^{d})$] the class of
$\reals^d$-valued functions on $E$ with continuous (resp., locally
$\gamma$-H\"{o}lder continuous) partial derivatives of $m$th order. The
superscripts are dropped for $m=0$ or $d=1$, so that
$C^{0,\gamma}(E,\reals^{1})$ is denoted by $C^{\gamma}(E,\reals)$.
The following assumption prescribes that the coefficients $r,\mu
,b,A,\Sigma$
and $\Upsilon$ are regular and nondegenerate.
\begin{assum}\label{regmod}
$r\in C^{\gamma}(E,\reals)$, $b\in C^{1,\gamma}(E,\reals^k)$,
$\mu\in C^{1,\gamma}(E,\reals^n)$, $A\in
C^{2,\gamma}(E,\reals^{k\times k})$, $\Sigma\in
C^{2,\gamma}(E,\reals^{n\times n})$ and $\ups\in
C^{2,\gamma}(E,\reals^{n\times k})$. The symmetric matrices $A$ and
$\Sigma$ are positive definite for all $y\in E$.
\end{assum}

The region $E$ is typically of the form $E = \reals^d\times
(0,\infty)^{k-d}$ for some $0\le d\le k$, as in the case of multivariate
Ornstein--Uhlenbeck processes, Feller diffusions or a combination
thereof. \citet{KarFern2005} consider models in which $E$ is the open simplex
in $\reals^{n-1}$.

To construct the processes $(R,Y)$, let $\Omega=
C([0,\infty),\reals^{n+k})$ endowed with the topology of uniform
convergence on bounded intervals. Let $\mathcal{B}$ be the Borel
sigma algebra and let $(\mathcal{B}_t)_{t\geq0}$ be the
filtration generated by the coordinate process $X$ defined by
$X_t(\omega) = \omega_t$ for $\omega\in\Omega$. For a second-order
differential operator $L$ of the form in \eqref{EPoperatordef}
below, a solution to the martingale problem for $L$ on
$\reals^n\times E$ is a family of Borel probability measures
$(P^{x})_{x\in\reals^n\times E}$ on $(\Omega,\mathcal{B})$ such
that, for each $x\in\reals^n\times E$, (i) $P^x(X_0 = x) = 1$, (ii)~%
$P^x(X_t\in\reals^n\times E, \forall t\geq0) = 1$ and (iii)~$(f(X_t)-f(X_0)-\int_0^t (L f)(X_u)\,du;\mathcal{B}_t)$ is a
$P^x$ martingale for all $f\in
C^2_0(\reals^n\times E)$.

The next assumption ensures that $\mu,b,A,\Sigma$ and $\Upsilon$
identify the
law of $(R,Y)$. $x=(r,y)$, with $r\in\reals^n, y\in E$, denotes
the starting points of $R$ and~$Y$.

\begin{assum}\label{assmod}
There exists a unique solution $(P^{(r,y)})_{r\in\reals
^n,y\in E}$ to the martingale problem for
%
%
\begin{eqnarray}\label{EPoperatordef}
L &=& \frac{1}{2}\sum_{i,j=1}^{n+k}\tilde{A}^{i,j}(x)\frac{\partial
^2}{\partial x_i\,\partial x_j} +
\sum_{i=1}^{n+k}\tilde{b}^{i}(x)\frac{\partial}{\partial x_i},
\nonumber
\\[-8pt]
\\[-8pt]
\nonumber
\tilde{A} &=& \pmatrix{
\Sigma& \Upsilon\vspace*{2pt}\cr
\Upsilon' & A},
\qquad
\tilde{b} = \pmatrix{
\mu\cr
b}.
\end{eqnarray}
\end{assum}

Since $R_0 = 0$ for all the models considered in this paper, the
measure $P^{(0,y)}$ in Assumption \ref{assmod} is simply denoted as
$P^y$. Henceforth, consider the filtration $(\mathcal{F}_t
)_{t\geq0}$ defined as $\mathcal{F}_t = \mathcal{B}_{t+}$, that is,
the right continuous envelope of $\mathcal{B}_{t}$. Under Assumption
\ref{regmod} for $f\in C^2_0(\mathbb R^n\times E)$,
$f(X_t) -f(X_0) - \int_0^t(Lf)(X_u)\,du$ is a martingale also with
respect to $\calf_t$; cf. \citet{MR1121940}, Section~5.4.

\begin{remark}
For consistency of notation, Assumption \ref{assmod} involves the joint
diffusion process $(R,Y)$. However, it is essentially an assumption on
the process $Y$. Indeed, if Assumption \ref{regmod} holds, and if there
is a unique solution $(Q^y)_{y\in E}$ to the martingale
problem for the operator
%
%
\begin{equation}\label{Eyoperatordef}
L^Y = \frac{1}{2}\sum_{i,j=1}^{k}A^{i,j}(y)\frac{\partial
^2}{\partial
y_i\,\partial y_j} +
\sum_{i=1}^{k}b^{i}(y)\frac{\partial}{\partial y_i}
,
\end{equation}
then Assumption \ref{assmod} holds. To see this fact, first consider
the martingale problem for
the operator associated to $(Y,B)$, where $B$ is an $n$-dimensional standard
Brownian motion independent of $Y$. Then, write $Z = \rho W + \bar
{\rho}B$
where $\bar{\rho}$ is a square root of $1_n-\rho\rho'$ and define the
integrals for $R$ in \eqref{reteq} accordingly.
\end{remark}

\subsection{Trading strategies}\label{sec1.2}
An investor trades in the market according to a~portfolio $\pi=
(\pi^i_t)^{1\leq i\leq n}_{t\ge0}$, representing the proportions of
wealth in each risky asset. Since the investor observes the state
variables $Y$ and the asset returns~$R$, the portfolio $\pi$ is
adapted to the filtration generated by $(R,Y)$, and is
$R$-integrable.\footnote{$R$-integrability means each of the integrals
$\int_0^{\cdot}\pi'\mu\,dt, \int_0^{\cdot}\pi'\sigma\,dZ_t$ is well
defined.} The corresponding wealth process $X^{\pi} =
(X^\pi_t)_{t\ge0}$ follows:
\begin{equation}\label{Estrategywealthdynamics}
\frac{dX^{\pi}_t}{X^{\pi}_t} = r(Y_t)\,dt + \pi_t'\,dR_t.
\end{equation}
%
Note that a positive initial capital $X_0\ge0$ implies a positive
wealth at all times, that is, $X^\pi_t\ge0$ a.s. for all $t\ge0$,
thereby ruling out doubling strategies [see, e.g.,
\citet{MR622165}].\vadjust{\goodbreak}

\subsection{Preferences}\label{sec1.3}
The investor's preferences are modeled by the power utility function
%
%
%
%
\begin{equation}\label{homo}
U(x) = \frac{x^p}{p},\qquad
{p<1, p\neq0}.
\end{equation}
Denoting by $E_P^y$ the expectation with respect to $P^y$, the goal
is to maximize expected utility from terminal wealth. With a finite
horizon $T$, the problem is
%
%
\begin{equation}\label{maxutfo}
\max_{\pi} \frac1pE_P^y[(X^{\pi}_T)^p].
\end{equation}
Since power utility is homothetic [$U(c x)=c^p U(x)$], it
suffices to consider the case $X_0=1$ of unit initial wealth.
Henceforth, $q$ denotes the conjugate exponent to $p$,
%
%
\begin{equation}\label{Eqdef}
q := \frac{p}{p-1}.
\end{equation}
To preserve economic intuition, recall that
$p=1-\gamma$ and $q=1-\frac1\gamma$, where $\gamma$ is the investor's
relative risk aversion.
Risk aversion increases as $q$ increases and as~$p$ decreases, and
logarithmic utility corresponds to the limit $p\rightarrow0$.

The martingale approach to utility maximization relies on the duality between
final payoffs and pricing rules, represented by the related concepts of
stochastic discount factors and martingale measures.
\begin{defin}
A stochastic discount factor is a strictly positive adapted process
$M=(M_t)_{t\ge0}$, such that $M S$ is a martingale
\[
E_{P}^y[M_t S^i_t | \calf_s]=M_s S^i_s \qquad\mbox{for all
}0\le
s\le t, 0\le i\le n.
\]
A martingale measure is a probability $Q$, such that $Q|_{\calf_t}$
and $P^y|_{\calf_t}$ are equivalent for all $t\in[0,\infty)$ and the
discounted prices $S^i/S^0$ $($or equivalently, the excess returns
$R^i)$ are $Q$-martingales for all $1\le i\le n$.
\end{defin}

Martingale measures and stochastic discount factors are in a
one-to-one correspondence through the relation
$\frac{dQ}{dP^y}|_{\calf_t}=S^0_t M_t$, although their
distinction is important in the present context of stochastic
interest rates. Except for a complete market in which the martingale
measure is unique, the description of a pricing rule requires the
choice of unhedgeable risk premia $\eta$.
For any $(\mathcal{F}_t)_{t\geq0}$ adapted,
integrable process $\eta$, a candidate (local) martingale measure
$Q^{\eta,y}$ is given by $dQ^{\eta,y}/dP^{y}|_{\mathcal{F}_t} =
Z^{\eta
}_t$, where
%
%
\begin{equation}\label{Epredictrepeqn}
Z^{\eta}_t =
\mathcal{E}\biggl(-\int_0^{\cdot}
(\mu'\Sigma^{-1}+\eta'\ups'\Sigma^{-1})\sigma\,dZ+ \int_0^{\cdot}
\eta' a \,dW\biggr)_t
\end{equation}
and where $\cale(X)_t=\exp(X_t-\frac12 \langle X\rangle_t)$. Clearly,
$Z$ must be a martingale for~$Q^{\eta,y}$ to be an equivalent local
martingale measure.\vadjust{\goodbreak} For such a $\eta$, let $M^{\eta}$ denote the
corresponding stochastic discount
factor
%
%
\begin{equation}\label{EMetadef}
 M^{\eta}_t = \exp\biggl(-\int_0^tr \,dt\biggr
)\mathcal{E}\biggl(-\int
_0^{\cdot}
(\mu'\Sigma^{-1}+\eta'\ups'\Sigma^{-1})\sigma\,dZ+ \int_0^{\cdot}
\eta' a \,dW\biggr)_t.\hspace*{-30pt}
\end{equation}
Note that for any strategy $\pi$ and risk premia $\eta$, by
\eqref{Estrategywealthdynamics}, the process $X^{\pi}M^{\eta}$ is a
super-martingale, even if the right-hand side of \eqref{Epredictrepeqn}
is only a local martingale. For power utility, the
following lemma applied to $X=X^{\pi}_T$
and $M = M^{\eta}_T$ for any $T>0$ shows that the duality bound is an
immediate consequence of H\"{o}lder's inequality and of the super-martingale
property of the process~$X^\pi M^\eta$.

\begin{lemma}\label{LHolderdualitybound}
Let $X,M$ be random variables on a probability space
$(\Omega,\mathcal{F},P)$ such that $X, M>0$ $P$-almost
surely and
$\espalt{P}{XM}\leq1$. Then
%
%
\begin{equation}\label{dualbound}
\frac{1}{p}\espalt{P}{X^p}\leq\frac{1}{p}\espalt{P}{M^q}^{1-p}
\end{equation}
and equality holds if and only if $\espalt{P}{X M} =1$ and, for some
$\alpha> 0$,
%
%
\begin{equation}\label{firstor}
X^{p-1}=\alpha M.
\end{equation}
\end{lemma}

Equation \eqref{dualbound} bounds the utility of any terminal wealth
by a moment of any stochastic discount factor and vice versa.
The first-order condition \eqref{firstor} is the usual alignment of
marginal utilities with state-price densities.

Consider\vspace*{-1pt} a finite horizon $T$. Lemma \ref{LHolderdualitybound}
implies that a pair $(\pit, \etat)$ of a portfolio~$\pit$ and risk
premia $\etat$ such that\vspace*{-2pt} $X=X^{\pit}_T$ and $M=M^{\etat}_T$ is
optimal if it satisfies \eqref{firstor} and $E_{P}^y[X^{\pit}_T
M^{\etat}_T]=1$. Denoting by $u_T(y)$ the value function, that is, the
maximal expected utility, the following equalities hold:
%
%
\begin{equation}\label{optrel}
\frac1pE_P^y[(X^{\pit}_T)^p]=u_T(y)=\frac1pE_P^y[(M^{\etat}_T)^q]^{1-p}
\end{equation}
hence, $\pit$ is the optimal portfolio and the stochastic discount
factor $M^{\etat}$ identifies the pricing rule that makes an
investor indifferent between buying and selling a small amount of
any payoff, including unhedgeable ones.

\subsection{Long-run optimality}\label{sec1.4}

In the Markov model defined by \eqref{rateeq}--\eqref{correq},
sto\-chastic control arguments [see, e.g., \citet{MR1922696}, among many others],
show that the pair $(\pit,\etat)$ achieving optimality is of the
form $\pi^T_t = \pit(T-t,Y_t)$ and $\eta^T_t = \etat(T-t,Y_t)$ for
deterministic functions
\[
\pit: [0,T]\times E\mapsto\reals^n\qquad\etat: [0,T]\times
E\mapsto\reals^k.
\]
Thus, optimal policies depend on both state variables and
the residual horizon. This joint dependence is the major source of
intractability in portfolio choice and derivatives pricing problems.

\citet{brandt1999epa}, \citet{barberis2000ilr} and \citet
{wachter2002pac} report that optimal policies converge rapidly to
functions of state variables alone.
Thus the optimal policy for a long horizon $[0,T]$ is approximately
equal to a time-homogeneous function for most of the interval,
departing from it as the horizon $T$ approaches. The question is
whether using a time-homogenous policy throughout the interval $[0,T]$
can be approximately optimal.

For any functions $\pi\in C(E,\reals^n)$, $\eta\in
C(E,\reals^k)$ consider the portfolio $\pi=
(\pi(Y_t))_{t\geq0}$ and risk premia $\eta= (\eta(Y_t))_{t\geq
0}$. At any finite horizon $T$, the duality bound \eqref{dualbound}
implies that
%
%
\begin{equation}\label{duaineq}
\frac1pE_P^y[(X^{\hatp}_T)^p]
\le u_T(y)\le
\frac1pE_P^y[(M^{\mfnct}_T)^q]^{1-p}.
\end{equation}
The first inequality reflects the potential gap between the utility of the
long-run portfolio and the value function. A tangible measure of this
gap is
the increase in the risk-free rate $l_T$ required to recover this loss,
as to match the expected utility of the long-run optimal portfolio
under the
higher rate with the maximum expected utility at the regular rate. This
is the
\textit{certainty equivalent loss}, defined as
%
%
\begin{equation}\label{cel}
\frac1pE_P^y[(e^{l_T T} X^{\hatp}_T)^p]=u_T(y).
\end{equation}
Substituting \eqref{cel} into \eqref{duaineq} yields an upper bound
on $l_T$
%
%
\begin{equation}\label{celbou}
l_T\le
\frac1p
\biggl(
{\frac1T\log E_P^y[(M^{\mfnct}_T)^q]^{1-p}-\frac1T\log
E_P^y[(X^{\hatp}_T)^p]}
\biggr).
\end{equation}
This argument motivates the definition of a pair $(\pi,\eta)$ as
long-run optimal when its certainty equivalent loss vanishes for
long horizons.
\begin{defin}
\label{lrdef}
A pair $(\pi,\eta)\in C(E,\reals^n)\times C(E,\reals^k)$ is
\textit{long-run optimal} if, for all $y\in E$,
%
%
\begin{equation}\label{Eprimalduallongtermoptimality}
\limsup_{T\conv\infty}\frac1p
\biggl({\frac1T\log E_P^y[(M^{\mfnct}_T)^q]^{1-p}-\frac1T\log
E_P^y[(X^{\hatp}_T)^p]}\biggr)=0.
\end{equation}
\end{defin}

Long-run optimality defined here is essentially equivalent to the
criterion used by \citet{grossman1992odt} to solve portfolio choice
problems with leverage constraints. \citet{grossman1993ois} apply
the same idea to drawdown constraints and \citet{dumas1991esd} to
transaction costs. The \textit{risk-sensitive control} literature
studies a similar objective for multidimensional linear diffusions.

Definition \ref{lrdef} departs from the existing literature by
examining both the primal (investment) and the dual (pricing)
problems. This perspective yields verification theorems that are
valid for general multidimensional diffusions, provides estimates on
finite-horizon performance and allows to identify the parameter
sets for which long-run optimality holds.

Definition \ref{lrdef} allows another interpretation based on
management fees: an investor with sufficiently long horizon prefers a
long-run optimal portfolio to the optimal finite-horizon portfolio
if the long-run portfolio has slightly lower fees. The argument is
straightforward: suppose that the portfolio $\hatp$ requires the
payment of a
(continuously compounded) fee $\varphi$, while the finite-horizon
portfolio $\hatp_T$ entails a higher\vspace*{1pt} fee $\varphi_T\ge
\varphi$.
Accounting for such fees, the portfolio~$\hatp$ has expected
utility $\frac1p E_P^y[(X^{\hatp}_T e^{-\varphi T})^p]$.
However, by the bound~\eqref{dualbound} the finite-horizon
portfolio~$\hatp_T$ satisfies
\[
\frac1pE_P^y[(X^{\hatp_T}_T e^{-\varphi_T T})^p]\le
\frac1pE_P^y[(M^{\mfnct}_T)^q]^{1-p}e^{-p\varphi_T T}.
\]
Hence, $\hatp$ is preferred to $\hatp_T$ when
%
%
\begin{equation}\label{cell}
\varphi_T-\varphi\ge
\frac1p
\biggl(
{\frac1T\log E_P^y[(M^{\mfnct}_T)^q]^{1-p}-\frac1T\log
E_P^y[(X^{\hatp}_T)^p]}
\biggr).
\end{equation}
For a long-run optimal pair $(\pi,\eta)$, the limit of the
right-hand side is zero. Thus, for any minimal difference in fees, a
long-run optimal portfolio is preferable for investors with
sufficiently long horizons.

\section{Long-run analysis}\label{seclongrun}

The construction of long-run optimal portfolios $\pi$ and risk premia
$\eta$
takes place in two steps. In the first step, Theorem \ref{Tpdemodel}
computes the finite-horizon performance of the long-run optimal
``candidates''~%
$\pi$ and $\eta$. In the second step, Theorem \ref{Tlongrunopt}
establishes a sufficient condition for long-run optimality, requiring
that the
bounds found in the first step converge at long horizons.

The candidate long-run optimal $\pi$ and $\eta$ crucially depend
upon the solution of the quasi-linear partial differential equation
(PDE) in \eqref{EmainvPDE}, which acts as a long-run version of the
Hamilton--Jacobi--Bellman equation. Thus, Theorems \ref{Tpdemodel}
and~\ref{Tlongrunopt} are akin to verification theorems of
stochastic control theory, but for the asymptotic objective in Definition
\ref{lrdef}.
An advantage of these results is that they only rely on the local properties
of the processes~$(R,Y)$, avoiding the knowledge of the transition
density of
$Y$, which may be very complicated if known at all.

The second part of this section studies the existence of solutions to the
ergodic Bellman equation in \eqref{EmainvPDE}. Theorems
\ref{TLambdastructure}, \ref{Tbellmansolution} and Proposition
\ref{Ptestfortightness} below adapt the results
of \citet{MR2206349} to the present setting and under some extra
conditions in
addition to the assumptions of Theorems \ref{Tpdemodel} and~%
\ref{Tlongrunopt}. Their main message is that the quasi-linear PDE
generally admits only one candidate for long-run optimality. As shown with
examples in Section~\ref{secapp}, this candidate may or may not be long-run
optimal.


\subsection{Main results}\label{sec2.1}

Recall that, although their dependence on $y$ is omitted to alleviate
notation, $b$, $\mu$, $\Sigma$, $\ups$, $A$, $\phi$ and $v$ are
functions of the state
variable~$y$.

\begin{theorem}\label{Tpdemodel}
In addition to Assumptions \ref{regmod} and \ref{assmod}, assume that:
\begin{longlist}[(ii)]
\item[(i)]\label{Itdiffexpression}
$v\in C^2(E,\reals)$ and $\lambda\in\reals$ solve the
ergodic HJB equation (cf. Section~\ref{sec3.2})
%
%
\begin{eqnarray}\label{EmainvPDE}
&&pr - \frac q2\mu'\Sigma^{-1}\mu+
\frac{1}{2}\nabla{v}'(A-q\ups'\Sigma^{-1}\ups)\nabla{v}
\nonumber
\\[-8pt]
\\[-8pt]
\nonumber
&&\qquad{} +
\nabla{v}'(b-q\ups'\Sigma^{-1}\mu) +
\frac{1}{2}\trace(AD^2v)= \lambda;
\end{eqnarray}
\item[(ii)]\label{Itwellposed}
there is a unique solution $(\hat{P}^{r,y})_{r\in\reals
^n,y\in E}$
to the martingale problem for
%
%
\begin{eqnarray}\label{EhatPoepratordef}
\hat{L} &=& \frac{1}{2}\sum_{i,j=1}^{n+k}\tilde{A}^{i,j}(x)\frac
{\partial
^2}{\partial x_i\,\partial x_j} +
\sum_{i=1}^{n+k}\hat{b}^{i}(x)\frac{\partial}{\partial x_i},
\nonumber
\\[-8pt]
\\[-8pt]
\nonumber
\hat{b} &= &\pmatrix{
\displaystyle\frac{1}{1-p}(\mu+\Upsilon\nabla v) \vspace*{2pt}\cr
b-q\Upsilon'\Sigma^{-1}\mu+
(A-q\Upsilon'\Sigma^{-1}\Upsilon)\nabla v
},
\end{eqnarray}
where $\tilde{A}$ is as in \eqref{EPoperatordef}.
\end{longlist}
Then, the pair $(\hatp,\mfnct)$ given by
\begin{equation}
\label{lrport} \hatp= \frac{1}{1-p}\Sigma^{-1}(\mu+
\ups\nabla v) ,\qquad
\mfnct= \nabla v
\end{equation}
satisfies the equalities
\begin{eqnarray}
\label{priest}
E_P^y[(X^{\hatp}_T)^p] &=& e^{\lambda T+v(y)}E_{\hat
P}^y\bigl[e^{-v(Y_T)}\bigr],\\
\label{duaest} E_P^y[(M^{\mfnct}_T)^q]^{1-p} &=& e^{\lambda
T+v(y)}E_{\hat P}^y\bigl[e^{-({1}/{(1-p)})v(Y_T)}\bigr]^{1-p}.
\end{eqnarray}

\end{theorem}

\begin{remark}
Equations \eqref{priest} and \eqref{duaest} provide lower and upper
bounds on finite-horizon expected utility. Indeed, the duality
inequality \eqref{dualbound} yields
\begin{eqnarray*}
\frac1p e^{\lambda T+v(y)}E_{\hat P}^y\bigl[e^{-v(Y_T)}\bigr]&= &
\frac1p
E_P^y[(X^{\hatp}_T)^p] \le u_T(y)\le \frac1p E_P^y[(M^{\mfnct}_T)^q]^{1-p}\\
&=& \frac1p e^{\lambda
T+v(y)}E_{\hat P}^y\bigl[e^{-({1}/{(1-p)})v(Y_T)}\bigr]^{1-p}.
\end{eqnarray*}
Combining \eqref{priest} and \eqref{duaest} with \eqref{celbou} yields
the central quantitative implication: an upper bound on the
certainty equivalent loss
%
%
\begin{equation}\label{cebou}
l_T\le\frac1p \biggl(
\frac1T\log{E_{\hat P}^y\bigl[e^{-({1}/{(1-p)})v(Y_T)}\bigr]^{1-p}}-
\frac1T\log{E_{\hat P}^y\bigl[e^{-v(Y_T)}\bigr]} \biggr).
\end{equation}

\end{remark}

Theorem \ref{Tpdemodel} now reduces the long-run optimality
(Definition \ref{lrdef}) of $(\hatp,\mfnct)$ to the condition that
the right-hand side in \eqref{cebou} converges to zero. Theorem
\ref{Tlongrunopt} below provides a criterion that covers most
applications and Proposition \ref{PfadsisOK} below shows a
model in which this criterion is sharp, in that it holds for all
the parameter values for which long-run optimality holds.

\begin{theorem}\label{Tlongrunopt}
If, in addition to the assumptions of Theorem \ref{Tpdemodel}:
\begin{longlist}[(ii)]
\item[(i)]
the random variables $(Y_t)_{t\geq0}$ are $\hat P^y$-tight\footnote
{Recall that a family of $E$-valued random variables $(X_t)_{t\geq0}$
is \textit{P-tight in $E$}\vspace*{-1pt} if the induced measures $(P\circ
X^{-1}_t)_{t\geq0}$ form a tight family in $M_1(E)$, the space
of Borel measures on $E$. Thus, $(X_t)_{t\geq0}$ is $P$-tight in $E$
if for each $\eps
> 0$ there exists a compact $K_{\eps}\subset E$ such that
$\sup_{t\geq0}P(X_t\in K^c_{\eps})\leq\eps$.}
in $E$ for each $y\in E$;

\item[(ii)]
$\sup_{y\in E} F(y)<+\infty$, where $F\in C(E,\reals)$ is defined
as
%
%
\begin{equation}\label{Elongrunoptub}
 F=
\cases{
\biggl(pr-\lambda-\displaystyle\frac{q}{2}\mu
'\Sigma^{-1}\mu+
\displaystyle\frac{q}{2}\nabla v'\ups'\Sigma^{-1}\ups\nabla
v\biggr)e^{-v},
\qquad p<0,\vspace*{2pt}\cr
\biggl(pr-\lambda-\displaystyle\frac{q}{2}\mu
'\Sigma^{-1}\mu
-\displaystyle\frac{q}{2}\nabla v'(A-\ups'\Sigma^{-1}\ups)\nabla
v\biggr)e^{-({1}/{(1-p)})v}, &\vspace*{2pt}\cr
\hspace*{237pt}0<p<1.}\hspace*{-35pt}
\end{equation}
\end{longlist}
Then the pair $(\hatp,\mfnct)$ in \eqref{lrport} is long-run
optimal.
\end{theorem}

Section \ref{secapp} shows how to check conditions (i) and (ii) in typical
classes of models.

\begin{remark} A sufficient condition for (i) above to hold is that there
exists a nonnegative $\psi\in C(E,\reals)$ such that, for each $n$,
the level set $K_n\equiv\{y\in E : \psi(y)\leq n\}$ is compact and that
$M\equiv\sup_{t\geq0} E_{\hat{P}}^y[\psi(Y_t)] < \infty$. If such
a $\psi$ exists, then Markov's inequality implies that, for each $n$
\[
\sup_{t\geq0}\hat{P}^y(Y_t\in K_n^c)\leq\frac{1}{n}\sup_{t\geq
0}E_{\hat{P}}^y[\psi(Y_t)] = \frac{M}{n}.
\]
Thus, $\hat{P}^y$ tightness in $E$ follows.
\end{remark}
%
\subsection{Solutions to the ergodic Bellman equation}\label{SSergodicbellman}
This subsection provides conditions for the existence of a solution pair
$(v,\lambda)$ to \eqref{EmainvPDE}, such that the tightness
condition in
Theorem \ref{Tlongrunopt} holds. These results are obtained adapting the
arguments in \citet{MR2206349} to the present setting.
Define $\Lambda$ as the set of $\lambda\in\reals$ for which a solution
$v$ to \eqref{EmainvPDE} exists,
%
%
\begin{equation}\label{Ebellmansolutionrange}
\Lambda= \{\lambda\in\reals\ |\ \exists v\in
C^{2,\gamma}(E,\reals) \mbox{ solving }
\eqref{EmainvPDE}\}.
\end{equation}
$\Lambda$ depends both on the region $E$ and on the coefficients in the
PDE \eqref{EmainvPDE}.
The foregoing results require the following assumption on the region $E$,
which holds in virtually all models in the literature.
\begin{assum}\label{regmoddomain}
There exist $y_0\in E$ and a sequence of bounded open subsets
$E_n\subset E$,
star-shaped\footnote{Recall that $F\subseteq\reals^k$ is \textit{star-shaped}
for some $x_0\in F$ if for each $x\in F$ the segment $\{\alpha
x_0 +
(1-\alpha)x;\ 0\leq\alpha\leq1\}$ is contained within $F$. A convex
set is star-shaped with respect to any of its points.} with respect to
$y_0$ and with a $C^{2,\gamma}$ boundary, and strictly increasing to $E$,
in that $E = \bigcup_{n=1}^{\infty} E_n$ and $\bar{E}_n\cap
(E\setminus
E_{n+1}) =\varnothing$.
\end{assum}

This assumption is satisfied by any convex set $E$ for which there is a~convex
function $\psi\in C^{2,\gamma}(E,\reals)$ such that $\psi
(y)\uparrow
\infty$ as
$y\conv\partial E$. In this case, it suffices to set $E_n = \{y\in E:
\psi(y)<
n\}$. The next assumption requires that the potential ${pr} -
\frac{q}{2}\mu'\Sigma^{-1}\mu$ is bounded from above. When it does
not hold,
solutions to \eqref{EmainvPDE} may not exist [\citet{MR1326606},
Chapter 4.5].
\begin{assum}\label{boundedpotential}
$\sup_{y\in E} (pr(y) - \frac{q}{2}\mu(y)'\Sigma
(y)^{-1}\mu
(y)) < \infty$.
\end{assum}

Note that Assumption \ref{boundedpotential} always holds if $p<0$, and the
interest rate~$r$ is bounded from below, which is a typical situation in
financial models.
With Assumptions~\ref{regmod}, \ref{assmod} and \ref{regmoddomain},
denote by $(R,Y)$ the
coordinate process of the solution $(P^y)_{y\in E}$ of the martingale problem
corresponding to the operator $L$ from \eqref{EPoperatordef}.
Regarding the
state variable $Y$, the statement of existence results requires a few
basic definitions in
ergodic theory [see \citet{MR1326606}, \citet{MR1152459} for more
details]. $Y$~is \textit{transient} if $P^y(Y$ is
eventually in $E\setminus E_n \mbox{ for all
}n\ge N(\omega)) =1$ for all $n\ge1$
and $y\in E$. $Y$ is \textit{recurrent} if $P^x(\tau(\eps,y)<\infty
)
= 1$ for all $x,y\in E$ and $\eps>0$, where $\tau(\eps,y) = \inf\{
t\geq
0\ | |Y_t - y|\leq\varepsilon\}$.
If $Y$ is recurrent, there exists some\vspace*{1pt} $\tilde{\phi} > 0$ such that
$\tilde{L}\tilde{\phi} = 0$, where $\tilde L$ is the formal adjoint
to $L$.
$Y$ is \textit{positive recurrent}, or \textit{ergodic}, if $\int
_E\tilde
{\phi}\,dy <
\infty$, and \textit{null recurrent} otherwise. If $Y$ is positive
recurrent and $\tilde{\phi}$ is normalized
to be a probability density, then for all $y\in E$ and $f\in
L^1(E,\tilde{\phi})$
%
%
\begin{equation}\label{Eergodicassumptions}
\lim_{T\uparrow\infty} E_{P}^y[f(Y_T)] =\int_E f \tilde{\phi}\,dx.
\end{equation}
If $Y$ is ergodic, then \eqref{Eergodicassumptions} implies that $Y$
is $P^y$-tight in $E$ for each $y\in E$, and for all $y\in E$ and
$(t_n)_{n\ge1}\uparrow\infty$, the measures $(P^y\circ Y_{t_n})^{-1}$
weakly converge to the measure with density $\tilde{\phi}$, which does
not depend upon the starting point $y\in E$.
With these definitions and results, the following theorem shows that
there exists only one possible pair $(\lambda,v)$, solving \eqref
{EmainvPDE}, which can lead to long-run optimality.

\begin{theorem}\label{TLambdastructure}
Let Assumptions \ref{regmod}, \ref{regmoddomain} and \ref
{boundedpotential}
hold. Then there exists $\lambda_c\in\reals$ such that $\Lambda=
[\lambda_c,\infty)$. Furthermore, $Y$ is $(\hat{P}^y)_{y\in
E}$-transient for any~$\hat P$ corresponding to a solution $(v,\lambda
)$ of \eqref{EmainvPDE} with $\lambda> \lambda_c$.
\end{theorem}

Clearly, $(\hat{P}^y)_{y\in E}$ transience and $(\hat{P}^y)_{y\in
E}$-tightness in $E$ are incompatible with one another.\vadjust{\goodbreak} Furthermore, since
ergodicity implies tightness,\vspace*{1pt} the question is whether the pair
$(\lambda_c,v_c)$ makes $Y$ ergodic under $(\hat{P}^y)_{y\in
E}$. The following results give conditions under which this is indeed
the case. The first proposition is valid for a single state variable and
constant correlations~$\rho$. In this case, \eqref{EmainvPDE} linearizes
under a power transformation, and classical tests for transience and
recurrence of one-dimensional diffusions apply. The second proposition
considers the general multidimensional case but under a~stronger restriction
on the potential term.

\subsubsection{One state, constant correlations}\label{SSonestate}

An important case leads to substantial simplifications.
\begin{assum}\label{onedelta}
Let Assumptions \ref{regmod} and \ref{boundedpotential} hold.
Further, assume that:

\begin{longlist}[(ii)]

\item[(i)]$E=(\alpha,\beta)$ with $-\infty\leq\alpha< \beta\leq
\infty$.

\item[(ii)]$\rho'\rho=\ups'\Sigma^{-1}\ups/A$ is constant.
\end{longlist}
\end{assum}

Note that Assumption \ref{regmoddomain} is always satisfied
for a single state variable.~Set
%
%
\begin{equation}\label{Edeltadef}
\delta= \frac{1}{1-q\rho'\rho}.
\end{equation}
The change of variable $\phi=\exp(v/\delta)$,
essentially equivalent to the power transformation of
\citet{MR1807876}, reduces the quasi-linear ODE in~\eqref{EmainvPDE} to the linear ODE
%
%
\begin{equation}\label{Ephislode}
\frac{A}{2}\dder{\phi}
+(b-q\ups'\Sigma^{-1}\mu)\der{\phi} +\frac1\delta
(V-\lambda)\phi=0.
\end{equation}

Let $\lambda\in\Lambda$ and let $\phi\in C^{2,\gamma}(E,\reals)$
with $\phi> 0$ be a solution to \eqref{Ephislode} obtained by $\phi
= \exp{(v/\delta)}$. Under
$\hat{P}^y$, $Y$ has the dynamics
%
%
\begin{equation}\label{EYhatpdynamics}
dY_t = \biggl(b - q\ups'\Sigma^{-1}\mu+
A\frac{\der{\phi}}{\phi}\biggr)\,dt + a\,d\hat W_t.
\end{equation}

Using Feller's test for explosions, the following proposition [\citet
{MR1326606}, Corollary 5.1.11] gives sufficient
conditions for $Y$ to be $\hat{P}^y$-tight in $E$ for the candidate optimal
pair $(\lambda_c,\phi_c)$.

\begin{proposition}\label{Ptestfortightness}
Let Assumption \ref{onedelta} hold, let $(\lambda_c,v_c)$ be as in
Theorem~\ref{TLambdastructure} and let $\phi_c = \exp(v_c/\delta)$. Denote by
%
%
\begin{equation}\label{Emnudef}
m_{\nu}(y) =
\frac{1}{A(y)}\exp\biggl(\int_{y_0}^y\frac{2(b-q\Upsilon'\Sigma
^{-1}\mu)(z)}{A(z)}\,dz\biggr),
\end{equation}
where $y_0\in(\alpha,\beta)$. Then, the family of random variables
$(Y_t)_{t\geq0}$ is $\hat{P}^y$-tight in $E$ if and only if
%
%
\begin{equation}\label{Eposreccrit}
\int_\alpha^{y_0} \frac{1}{\phi_c^2 A m_{\nu}} \,dy =
\int_{y_0}^\beta\frac{1}{\phi_c^2A m_{\nu}} \,dy = \infty
\quad\mbox{and}\quad
\int_\alpha^\beta\phi_c^2 m_{\nu}\,dy < \infty.
\end{equation}
%
\end{proposition}

\subsubsection{The general case}

This subsection treats the general case of $k$ state variables under
the following assumption.
\begin{assum}\label{generalcase}
There exists a function $w\in C^{2,\gamma}(E,\reals)$ such that
%
%
\begin{eqnarray}\label{Ewdropoff}
&&\lim_{n\uparrow\infty}\sup_{E\setminus E_n}\biggl({pr}-\frac{q}{2}\mu
'\Sigma^{-1}\mu
+ \frac{1}{2}\nabla{w}'(A-q\ups'\Sigma^{-1}\ups)\nabla{w}
\nonumber
\\[-8pt]
\\[-8pt]
\nonumber
&&\hspace*{82pt}{}+ \nabla{w}'(b-q\ups'\Sigma^{-1}\mu) +
\frac{1}{2}\trace(AD^2w)\biggr) = -\infty.
\end{eqnarray}

\end{assum}

\begin{remark}
For general regions $E$, condition \eqref{Ewdropoff} plays a similar
role as
condition (A3) in \citet{MR2206349}. Assumption \ref{generalcase} is
satisfied, for example, when
%
%
\begin{equation}\label{Ewdropoffalt}
\lim_{n\uparrow\infty}\sup_{E\setminus E_n}\biggl(
{pr} -
\frac{q}{2}\mu
'\Sigma^{-1}\mu\biggr) = -\infty.
\end{equation}
In this case, $w\equiv1$ satisfies \eqref{Ewdropoff}.
\end{remark}

In this setting, the main existence criterion is the following.
\begin{theorem}\label{Tbellmansolution}
Let Assumptions \ref{regmod}, \ref{regmoddomain}, \ref{boundedpotential}
and \ref{generalcase} hold and let $(\lambda_c,v_c)$ be as in Theorem
\ref{TLambdastructure}. Then $v_c$ is unique up to an additive
constant and
$(Y_t)_{t\geq0}$ is $\hat{P}^y$-tight in $E$ for all $y\in E$.
\end{theorem}

Proposition \ref{PfadsisOK} in Section \ref{secapp} below shows that
long-run optimality may still fail, even when the tightness condition
is satisfied
for the pair $(\lambda_c,v_c)$. The reason is that, even if $Y$ is ergodic
under $(\hat{P}^y)_{y\in E}$ with invariant density~$\tilde{\phi}$, the
ergodic property in \eqref{Eergodicassumptions} may not hold for
\eqref{priest} and \eqref{duaest} because the functions
%
%
\begin{equation}\label{tempor}
\exp(-v(y)),\qquad\exp\bigl(-(1-p)^{-1}v(y)\bigr)
\end{equation}
therein may not be in $L^1(E,\tilde{\phi})$. Thus, long-run optimality
requires additional assumptions, such as \eqref{Elongrunoptub} in
Theorem \ref{Tlongrunopt}.

If the functions in \eqref{tempor} are in $L^1(E,\tilde{\phi})$, then
\eqref{Eergodicassumptions} yields additional information about the
speed at
which the certainty equivalent loss $l_T$ converges to zero in the
limit of a
long horizon. The following proposition provides such a result, in the
case of
a single state variable. Recall that risk aversion is $\gamma= 1-p =
1/(1-q)$. The main message is that long-run optimality may fail
only for (i) high risk aversion and highly incomplete
market or
(ii) low risk aversion and nearly complete market. In
particular, for
$-1<q\le1/2$ (i.e., risk aversion within $1/2$ and $2)$, long-run optimality
holds for any level of incompleteness. In addition, for $q>1/2$
(resp., $q<-1)$, long-run optimality holds if $\rho'\rho<1/2$
(resp., $\rho'\rho>1/2)$.

Section \ref{secapp} takes up this issue in specific models, obtaining
necessary and sufficient conditions. By contrast, the following sufficient
condition holds under general assumptions, regardless of the model
considered.

\begin{proposition}\label{Ponestatenice}
Let Assumption \ref{onedelta} hold, and assume that $m_{\nu}$ in
\eqref{Emnudef} satisfies $\int_Em_{\nu}\,dy < \infty$. If
$(\lambda
_c,v_c)$ are
such that $Y$ is $(\hat{P}^y)_{y\in E}$-ergodic, then long-run
optimality holds if
%
%
\begin{equation}\label{Eqrhocondnice}
\rho'\rho\in
\cases{
\biggl[0,\displaystyle\frac{1}{2q}\biggr], & \quad$\mbox{for }\displaystyle\frac12 < q <
1,$\cr
[0,1] , & \quad$\mbox{for } -\!1 \leq q\le \displaystyle\frac12,q\neq0,$\cr
\biggl[\displaystyle\frac{1+q}{2q},1\biggr], &\quad$\mbox{for }q<-1$.}
\end{equation}
In such cases, there exists a constant $K>0$ such that the certainty
equivalent loss~$l_T$ satisfies
%
%
\begin{equation}\label{ECLEdecaynice}
0\leq\limsup_{T\uparrow\infty} T l_T \leq K.
\end{equation}
\end{proposition}

The next corollary states in an important special case, which does
not even require the knowledge of the principal eigenfunction $\phi$,
since $m_\nu$ only depends on the model parameters.
\begin{corollary}
Under Assumption \ref{onedelta}, if $\int_E m_{\nu}\,dy < \infty$ and
\eqref{Ewdropoffalt} are satisfied, then long-run optimality holds
for $q$
and $\rho'\rho$ satisfying \eqref{Eqrhocondnice}.
\end{corollary}

\section{Implications and ramifications}\label{interpret}

\subsection{The myopic probability}\label{sec3.1}
The bounds \eqref{priest} and \eqref{duaest} in Theorem
\ref{Tpdemodel} and assumption (i) in Theorem
\ref{Tlongrunopt} depend on the equivalent
probability $\hat P^y$, which plays a pivotal role in long-run
analysis. In general, $\hat P^y$ is neither the physical
probability~$P^y$ nor
a risk-neutral probability. Instead,\vspace*{1pt} its interpretation becomes
clear from its dynamics, which is (for $\hat P$-Brownian motions~$\hat
Z, \hat W$)
\begin{eqnarray*}
dR_t &=& \frac{1}{1-p}(\mu+ \ups\nabla v)\,dt
+\sigma\, d\hat Z_t,\\
dY_t &=& \bigl(b - q\ups'\Sigma^{-1}\mu+
(A-q\ups'\Sigma^{-1}\ups)\nabla v\bigr)\,dt + a\,d\hat W_t.
\end{eqnarray*}
Compare the original model, with price dynamics under $P^y$ and
power\vspace*{1pt} utility~${x^p}/{p}$, to the auxiliary model under $\hat P^y$
with logarithmic utility. The long-run\vadjust{\goodbreak} optimal portfolio in the two
models coincide. The first one is simply in \eqref{lrport}, while
the second one follows from the usual formula\vspace*{-3pt} $\pi=\Sigma^{-1}\hat
\mu$, where $\hat\mu=\frac{1}{1-p}(\mu+ \ups\nabla v)$
are the expected returns under $\hat P^y$. Thus a~long horizon,
power-utility investor under the probability $P^y$ behaves exactly
as a myopic (or logarithmic) investor under~$\hat P^y$.

This observation shows that $\hat P^y$ corresponds to the long horizon
limit of the probability $\mathbf R$ considered by \citeauthor
{MR2260066} (\citeyear{MR2260066,MR2288717,MR2438000}) in finite
horizon, in the context of sensitivity analysis pricing of option prices.
\citet{MR2330978} study mean--variance hedging for semimartingales
and obtain optimal strategies in terms of the predictable
characteristics of asset prices\vspace*{1pt} under an \textit{opportunity neutral}
probability $P^*$, which is similar in spirit to $\hat P^y$, in that
it reduces the mean--variance objective to a logarithmic utility
objective.

\subsection{Connections with stochastic control}\label{sec3.2}

Since the work of \citet{merton1969lps}, most of the dynamic portfolio choice
literature has employed stochastic optimal control as its main
analytical tool. The relation between Theorems~\ref{Tpdemodel},~%
\ref{Tlongrunopt} and the stochastic control approach becomes clear
by comparing equation~\eqref{EmainvPDE} to the Hamilton--Jacobi--Bellman (HJB) equations
of the utility maximization problem~\eqref{maxutfo}. Its value
function $u(x,y,t)$ depends on the current wealth $x$, the current
state $y$ and time $t$. The homogeneity of power utility entails
that $u(x,y,t)=x^p e^{w(y,t)}/p$, thereby removing wealth from the
reduced value function $w$. The corresponding HJB equation becomes
[see, e.g., \citet{MR1922696}]
%
%
\begin{eqnarray}\label{EmainHJBequation}
-\frac{\partial w}{\partial t}&=&pr - \frac{q}{2}\mu
'\Sigma^{-1}\mu+
\frac{1}{2}{\nabla{w}'}(A-q\ups'\Sigma^{-1}\ups){\nabla{w}}
\nonumber
\\[-8pt]
\\[-8pt]
\nonumber
&&{}+ {\nabla{w}'}(b-q\ups'\Sigma^{-1}\mu) +
\frac{1}{2}\trace(AD^2w)
\end{eqnarray}
with the terminal condition $w(y,T)=0$. Instead, the main PDE
\eqref{EmainvPDE} is
\begin{eqnarray*}
\lambda&= &pr - \frac{q}{2}\mu'\Sigma^{-1}\mu+
\frac{1}{2}\nabla{v}'(A-q\ups'\Sigma^{-1}\ups)\nabla{v}\\
&&{} +
\nabla{v}'(b-q\ups'\Sigma^{-1}\mu) +
\frac{1}{2}\trace(AD^2v).
\end{eqnarray*}
In the former equation, the unknown function $w$ depends on both
time $t$ and the state $y$, while $v$ in the latter equation only
depends on the state, although the constant $\lambda$ is also
unknown. Indeed, the former equation reduces to the latter under the restriction
\[
w(t,y)=\lambda(T-t)+v(y).
\]
%
This restriction gains analytical tractability by reducing the
dimension of the problem. The price of the tractability gain is that
solutions of the time-homogeneous equation in general do not satisfy
the boundary condition and, therefore, are not exactly optimal at any
time-horizon (except in the trivial case $v=0$, arising with
logarithmic utility or constant investment opportunities).

A special case of equation \eqref{EmainvPDE} appears in the
risk-sensitive control approach to optimal investment initiated by
\citet{MR1675114}. In a~linear diffusion model, they study the
problem
%
%
\begin{equation}\label{rsc}
\max_{\pi}\liminf_{T\rightarrow\infty}\frac1T\log E[(X^\pi_T)^p],
\end{equation}
where the supremum is taken over all progressively measurable
strategies. Risk-sensitive control relies on control techniques to
establish the existence and uniqueness to the homogeneous equation,
then attempts to establish its optimality in the sense of
\eqref{rsc}. \citeauthor{MR1802598}
(\citeyear{MR1802598,MR1910647}) carry out this program under the assumption
that $|p|$ is small, that is, if risk aversion is close enough to the
logarithmic case. The results in this paper, which apply to general nonlinear
models, shed new light on this literature by characterizing finite-horizon
performance. For example, Proposition \ref{PfadsisOK} below relaxes the
restriction of $|p|$ small to a necessary and sufficient condition and
explains the economic intuition behind it.

\subsection{Long-run $q$-optimal measure}\label{sec3.3}

For each value of the risk-aversion parameter $1-p$, the risk premia
$\mfnct$ in \eqref{lrport} deliver a pricing rule for derivatives
involving the partially unhedgeable state variable $Y$.
The martingale measure $Q^{\mfnct}$ corresponding to the risk premia
$\mfnct$ is a long-run version of the \textit{minimax martingale
measure} of \citet{MR1122311}, called $q$-optimal measure by
\citet{MR2092922} and \citet{MR2116796}. Its
formal dynamics is
%
%
\begin{equation}\label{qopt}
\cases{
dR_t =
\sigma\, d\tilde Z_t,\vspace*{2pt}\cr
dY_t = \bigl(b - \ups'\Sigma^{-1}\mu+ (A -
\ups'\Sigma^{-1}\ups)\nabla v\bigr)\,dt + a\,d\tilde W_t}
\end{equation}
for some Brownian motions $\tilde Z$ and $\tilde W$. Since this
dynamics is distinct from the one under $P^y$ and $\hat P^y$, in
general it is necessary to check its well-posedness, in the form of
Assumption \ref{assmod}.

Observe that the drift of $Y$ under the $q$-optimal measure has
three components. The first term $b$ is the drift under the original
measure $P^y$. The second term $\ups'\Sigma^{-1}\mu$ is the
risk-neutral adjustment due to the correlation between the returns
and the state shocks.
The last term $(A - \ups'\Sigma^{-1}\ups)\nabla v$
accounts for preferences which enter the equation through $v$.

\subsection{Complete and fully incomplete as duals}\label{sec3.4}

The formulas in \eqref{lrport} highlight the symmetric aspects of
complete markets, where $A=\ups'\Sigma^{-1}\ups$ identically and fully
incomplete markets, where $\ups=0$. In a complete market the pricing
problem is trivial, as the dynamics in \eqref{qopt} becomes
independent of the preference parameter~$p$,
%
%
\begin{equation}\label{qcom}
\cases{
dR_t =
\sigma\, d\tilde Z_t,\vspace*{2pt}\cr
dY_t = (b - \ups'\Sigma^{-1}\mu)\,dt + a\,d\tilde
W_t.}
\end{equation}
However, the investment problem is nontrivial because the optimal
portfolio includes a component that perfectly hedges the state
variables.

Conversely, in a fully incomplete market, myopic portfolios are
optimal because all portfolios evolve orthogonally to state
variables. However, a latent hedging motive remains present and
generates nonzero risk premia for state variables, which have a
potential as hedging instruments.
In both cases, it is market dynamics, and not preferences, which
make either the investment or the pricing problem trivial. By
contrast, logarithmic preferences or constant $\mu$ and $\sigma$
remove the intertemporal hedging motive entirely, making both
problems trivial.

The quasi-linear ODE in \eqref{EmainvPDE} becomes linear under a
transformation in both the complete and fully incomplete cases,
as in the one state variable case discussed in Section \ref
{SSonestate}. In the
complete case, the transformation $\phi= e^{(1-q)v}$ leads to
the linear equation
%
%
\begin{equation}\label{Elinearcomplete}
\hspace*{28pt}\frac{1}{2}\trace(AD^2\phi) +
\nabla\phi'(b-q\ups'\Sigma^{-1}\mu) +
(1-q)\biggl(pr-\frac{q}{2}\mu'\Sigma^{-1}\mu
-\lambda
\biggr)\phi= 0.
\end{equation}
In the fully incomplete case, the transformation $\phi= e^v$
leads to the linear equation
%
%
\begin{equation}\label{Elinearfullyincomplete}
\frac{1}{2}\trace(AD^2\phi) + \nabla\phi'b +
\biggl(pr-\frac{q}{2}\mu'\Sigma^{-1}\mu-\lambda
\biggr)\phi= 0.
\end{equation}
The criticality theory of \citet{MR1326606} applies to these cases
under multivariate restrictions similar to those given in Assumption
\ref{onedelta}, with $\delta= \frac{1}{1-q}$ in the complete case
and $\delta= 1$ in the fully incomplete case. Furthermore, the
multivariate results of Theorem \ref{Tbellmansolution} apply
as long as \eqref{Ewdropoff} can be verified.

\subsection{Long-run decomposition}\label{sec3.5}

The bounds \eqref{priest} and \eqref{duaest} decompose expected
utility and its dual into a common ``long-run'' component
$e^{\lambda T}$, and two ``transient'' components, in a close
analogy to \citet{hansenltr}. For a multiplicative functional $N$
of a Markov process $Y$, they propose the decomposition
%
%
\begin{equation}
\label{hsd} N_t=\exp(\rho t)\frac{\varphi(y)}{\varphi(Y_t)}\hat N_t,
\end{equation}
where $\rho$ and $\varphi$ are, respectively, the principal eigenvalue
and eigenfunction of the infinitesimal generator of $Y$ and $\hat
N_t$ is a martingale.
The bounds~\eqref{priest} and~\eqref{duaest} yield similar
expressions for terminal utilities and their dual counterparts
\[
{(X^{\hatp}_T)^p} = e^{\lambda T} \frac{e^{v(y)}}{e^{v(Y_T)}}
\frac{d\hat P^y}{d P^y} \qquad{(M^{\mfnct}_T)^q} =
\biggl(e^{\lambda
T}\frac{e^{v(y)}}{e^{v(Y_T)}}\biggr)^{{1}/{(1-p)}} \frac{d\hat
P^y}{d P^y}.
\]

These decompositions are precisely of the form in \eqref{hsd}, with
the difference that on the dual side the transient components
are powers of $e^v$, as opposed to $e^v$ itself.
Note also that the operator in \eqref{EmainvPDE} is not the
generator of $Y$ under either $P^y$ or~$\hat P^y$, since it is
nonlinear, has nonzero potential and has a different drift.
Further, the interpretation of $e^{\lambda T}$ as a long-run component hinges
on the condition that the $\hat P^y$-expectation of transient
components has a less than exponential growth, which means that
long-run optimality holds. This is not always the case; the examples
in Section \ref{secapp} show how parameter restrictions are
necessary even in the most common models.

\subsection{Large deviations}\label{sec3.6}

Theorems \ref{Tpdemodel} and \ref{Tlongrunopt} are closely
related to the results of \citeauthor{MR0386024}
(\citeyear{MR0386024,MR0428471,MR690656}) on large deviations of
occupation times for diffusions. Though the results also hold in the
multidimensional case of $k > 1$ state variables, the following
discussion considers a~single state variable for simplicity of
notation.

Let $E=(\alpha,\beta)$ for $-\infty\leq\alpha< \beta\leq\infty$
and consider a diffusion $Y$ with generator $L^Y$ from
\eqref{Eyoperatordef} (with $k=1$), assuming that the coefficients
$A$ and $b$ are such that $Y$ is positive recurrent under $(P^y)_{y\in
E}$. Let $m$ be the invariant measure which has a density by
\eqref{Eergodicassumptions}. With a slight abuse of notation, let
$m(dy) = m(y)\,dy$.

Denote by $M_1(\alpha,\beta)$ the space of Borel probability
measures on $(\alpha,\beta)$. Under certain conditions
on $Y$, \citeauthor{MR0386024} show that, for all continuous bounded
functions $V:(\alpha,\beta)\mapsto\reals$ and all
$y\in(\alpha,\beta)$,
%
%
\begin{equation}\label{dvas}
\qquad\lim_{T\rightarrow\infty} \frac1T \log E_P^y\biggl[\exp\biggl
(\int_0^T
V(Y_t)\,dt\biggr)\biggr]= \sup_{\mu\in
M_1((\alpha,\beta))}\biggl(\int_\alpha^\beta Vd\mu-I(\mu)\biggr).
\end{equation}
For $\mu\in M_1(\alpha,\beta)$ absolutely continuous with
respect to $m$\vspace*{1pt} (and hence, the Lebesgue measure) and with density
$\mu(y)$ such that $\psi(y)^2 \equiv\mu(y)/m(y)$ satisfies certain
regularity and decay conditions, the \textit{rate function} $I(\mu)$
reduces~to
\[
I(\mu)= \frac{1}{2}\int_\alpha^\beta A(\der{\psi})^2m\,dy.
\]
Using this representation for the rate function, the following
heuristic argument shows the relation between the \citet{MR0386024}
theory and long-run optimality. Consider the terminal utility of a
portfolio $\pi(Y_t)$ for some function $\pi: E\mapsto\reals^n$,
\[
(X^{\pi}_T)^p = \exp\biggl(\int_0^T\biggl(pr +
p\pi'\mu+ \frac
{1}{2}p(p-1)\pi
'\Sigma\pi\biggr)\,dt\biggr)
\mathcal{E}\biggl(\int_0^{\cdot}p \pi'\sigma\,dZ_t\biggr)_T.
\]
Define $P^y_\pi$ by setting $\frac{dP^y_\pi}{dP^y}$ equal to the
stochastic exponential in the last term of this equation. It follows
that
\[
E_P^y[(X^{\pi}_T)^p]= E_{P_\pi}^y\biggl[\exp\biggl(\int_0^T\biggl
(pr
+ p\pi'\mu+ \frac{1}{2}p(p-1)\pi'\Sigma\pi\biggr)\,dt\biggr
)\biggr].
\]
Assuming they may be applied under $P^y_{\pi}$, the
\citeauthor{MR0386024} asymptotics~\eqref{dvas} yield
%
%
\begin{eqnarray}\label{intob}
\quad&&\lim_{T\rightarrow\infty} \frac1T \log E_P^y[(X^{\pi}_T)^p]
\nonumber\\
&&\qquad= \sup_{\psi\in L^2_1(m_{\pi})} \int_\alpha^\beta\biggl
(\biggl({pr} +
p\pi'\mu+ \frac{1}{2}p(p-1)\pi'\Sigma\pi\biggr)\psi^2\\
&&\hspace*{179pt}\qquad\quad{}-\frac{1}{2}A(\der{\psi})^2\biggr
)m_{\pi}\,dy,\nonumber
\end{eqnarray}
where $m_\pi$ is the invariant density of $Y$ under $P^y_\pi$ and
$L^2_1(m_{\pi})$ is the unit disc in $L^2(m_{\pi})$. In a
similar manner to \eqref{Emnudef} in Proposition \ref{Ptestfortightness},
$m_{\pi}$ admits the formula
\[
m_{\pi}(y) =
\frac{1}{A(y)}\exp\biggl(\int_{y_0}^{y}\frac{2(b+p\Upsilon'\pi
)(z)}{A(z)}\,dz\biggr).
\]
Here, $y_0$ is some interior point in $E$. To make the
dependence on the portfolio~$\pi$ explicit, the change of
variable $\psi^2 m_\pi= \phi^2 m_{\nu}$ yields [$m_{\nu}$ is
defined in
\eqref{Emnudef}]
%
%
\begin{equation}\label{subst}
\frac{\der{\psi}}{\psi} = \frac{\der{\phi}}{\phi} +
\frac{\der{m}_{\nu}}{2m_{\nu}} - \frac{\der{m}_\pi}{2m_\pi} =
\frac{\der{\phi}}{\phi} - A^{-1}(q\ups'\Sigma^{-1}\mu+ p
\ups'\pi).
\end{equation}
Substituting \eqref{subst} into \eqref{intob}, the utility growth
rate becomes
%
%
\begin{equation}\label{clumu}
\sup_{\phi\in L^2_1(m_{\nu})}\int_\alpha^\beta(\pi'\mathbf
{A}\pi
+ \pi'\mathbf{B} + \mathbf{C})\phi^2m_{\nu}\,dy,
\end{equation}
where
\begin{eqnarray*}
\mathbf{A} &=& \frac{1}{2}p(p-1)(\Sigma- q\ups
A^{-1}\ups'), \qquad\mathbf{B} = p\biggl((1_{n}-q\ups
A^{-1}\ups'\Sigma^{-1})\mu+
\ups\frac{\der{\phi}}{\phi}\biggr),\\
\mathbf{C} &=& pr - \frac{1}{2}A\biggl(\frac{\der
{\phi
}}{\phi}\biggr)^2
+ q\mu'\Sigma^{-1}\ups\frac{\der{\phi}}{\phi} -
\frac{1}{2}q^2\mu'\Sigma^{-1}\ups A^{-1}\ups'\Sigma^{-1}\mu.
\end{eqnarray*}
The integrand is a quadratic function of $\pi$ and achieves its
optimum at
%
%
\begin{equation}\label{cannpi}
\hatp=\frac{1}{1-p}\Sigma^{-1}\biggl(\mu+
\delta\ups\frac{\der{\phi}}{\phi}\biggr).
\end{equation}
Thus, substituting \eqref{cannpi} into \eqref{clumu}, the utility
growth rate reduces to
%
%
\begin{equation}\label{Eoptimalpibound}
\sup_{\phi\in L^2_1(m_{\nu})}\int_\alpha^\beta\biggl(\biggl
({pr} -
\frac q2\mu'\Sigma^{-1}\mu\biggr)\phi^2 -
\frac{\delta}{2}A(\der{\phi})^2\biggr)m_{\nu}\,dy.
\end{equation}
A similar reasoning on stochastic discount factors delivers the
candidate long-run risk premia. The Euler--Lagrange equation
associated to \eqref{Eoptimalpibound} is the ODE in
\eqref{Ephislode}. Thus, large deviations arguments act as
a guide for producing the candidate long-run optimal policies.

This argument, which explains the formal connection with large
deviations, is suggestive but only heuristic. The main reason is
that the Donsker--Varadhan asymptotics are correct under some
delicate conditions which may fail to hold even in the simplest
models.

\section{Applications}\label{secapp}

This section applies the main results to two models, assuming that
the investor is more risk averse than the log investor (\mbox{$p<0$}).
In the first model, the state variables follow a multivariate
Ornstein--Uhlenbeck process that drives the drift of the return
process. Under general conditions, this model admits a unique
solution $v$ to \eqref{EmainvPDE}, leading to $\hat{P}^y$
tightness. For a~single state variable, long-run optimality is
characterized in terms of precise parameter restrictions.

In the second model, interest rates, drifts and volatilities are
stochastic. Each of these quantities is affine in a single common
state variable, which follows a Feller diffusion. Although this
model does not belong to the affine class, the long-run optimal
portfolios and risk premia have very simple expressions.

The parametric restrictions required by long-run optimality lead
under each single variate model to the same economic interpretation.
Long-run optimality does not hold at the conjunction of three
extreme situations: (i)~high covariation between risk premia and state
variables, (ii) nearly complete markets {and} (iii) high
risk-aversion. To understand this phenomenon, recall that long-run
optimality means that a time-homogenous strategy is approximately
optimal on a long time interval. Thus, the sub-optimality of the
long-run strategy in the latest part of the interval must lead to a
small utility loss. Since the myopic component of the optimal
finite-horizon portfolio is time-homogenous, any loss in utility is
attributed to the intertemporal hedging component. All of the three extreme
situations mentioned above concur to amplify the intertemporal hedging
component. First, the covariation of risk
premia is proportional to the hedging portfolios $\Sigma^{-1}\ups$.
Second, intertemporal hedging is more attractive in a nearly
complete market, where state variables are almost replicable. Third,
intertemporal hedging is higher for more risk-averse
investors who reduce long-term risk at the expense of short-term
return.

\subsection{Linear diffusion}\label{sec4.1}

This is the most common multivariate model, with
constant covariance matrices $\Sigma, \ups, A$ and drifts $r, \mu, b$
that are affine functions of the state variable. The dynamics is
%
%
\begin{equation}\label{lindif}
\cases{
dR_t = (\mu_0+\mu_1 Y_t)\,dt + \sigma\,dZ_t,\vspace*{2pt}\cr
dY_t = -bY_t \,dt+a\,dW_t,\vspace*{2pt}\cr
d\langle R,Y\rangle_t = \rho\,dt,\vspace*{2pt}\cr
r(Y_t) = r_0+r_1'Y_t,}
\end{equation}
where $\mu_0 \in\reals^n$, $\mu_1\in\reals^{n\times k}$, $\sigma
\in
\reals^{n\times n}$, $b\in\reals^{k\times k}$, $a\in\reals
^{k\times
k}$, $\rho\in\reals^{n\times k}$, $r_0\in\reals$ and $r_1\in
\reals
^k$. Under this model, $E=\reals^k$ and state variables follow a
multivariate Ornstein--Uhlenbeck
process. This setting is considered in most of the literature in risk
sensitive control mentioned in the \hyperref[intr]{Introduction} and
it is also implicit
in the use of vector autoregressions in the econometrics literature.
The coefficients in~\eqref{lindif} satisfy the following.
\begin{assum}\label{assmodlindif}
All four matrices $\Sigma=\sigma\sigma'$, $\mu_1'\mu_1$, $b + b'$ and
$a$ are positive definite.
\end{assum}

Assumption \ref{assmodlindif} implies that condition \eqref{Ewdropoffalt}
holds for $p<0$. Hence, by Theorem~\ref{Tbellmansolution} [or by the
results in \citet{MR2206349}], there is a unique pair $(\lambda,v)$
solving \eqref{EmainvPDE} such that for each $y\in E$, $(Y_t)_{t\geq
0}$ is $\hat{P}^y$-tight in $E$. The next theorem shows that $v$ is in
fact quadratic.

\begin{theorem}\label{teolin}
Let Assumption \ref{assmodlindif} hold for the model in \eqref
{lindif}. For
$p<0,$ equation \eqref{EmainvPDE} admits a unique solution $v$ such
that for
each $y\in\reals^k$, $(Y_t)_{t\geq0}$ is $\hat{P}^y$-tight in
$\reals^k$. The solution is of the form $v(y)=v_0'y - \frac
{1}{2}y'v_1 y$,
where the symmetric matrix $v_1\in\reals^{k\times
k}$ and the vector $v_0\in\reals^k$
and satisfy the algebraic equations
\begin{eqnarray}
\label{vzeq}
&&\bigl(v_1(A-q\ups'\Sigma^{-1}\ups)+(b+q\ups'\Sigma^{-1}\mu_1)'\bigr)v_0-pr_1
\nonumber
\\[-8pt]
\\[-8pt]
\nonumber
&&\qquad{}+q(\mu_1'-v_1\ups')\Sigma^{-1}\mu
_0=0,\\
\label{vueq}&&v_1(A-q\ups'\Sigma^{-1}\ups)v_1+v_1(b+q\ups'\Sigma^{-1}\mu_1)+
(b+q\ups'\Sigma^{-1}\mu_1)'v_1
\nonumber
\\[-8pt]
\\[-8pt]
\nonumber
&&\qquad{}-q\mu_1'\Sigma^{-1}\mu_1=0.
\end{eqnarray}
The corresponding utility growth rate equals
%
%
\begin{eqnarray}\label{vleq}
\lambda&=&pr_0-\frac{q}{2}\mu_0'\Sigma^{-1}\mu
_0+\frac
{1}{2}v_0'(A-q\ups
'\Sigma^{-1}\ups)v_0
\nonumber
\\[-8pt]
\\[-8pt]
\nonumber
&&{}-q v_0'\ups'\Sigma^{-1}\mu_0-\frac{1}{2}\trace(A v_1).
\end{eqnarray}
\end{theorem}

Equation \eqref{vueq} is a quadratic equation in the unknown matrix
$v_1$ and it belongs to the class of matrix Riccati equations
which arise in filtering theory and dynamical\vadjust{\goodbreak} systems. It does not
admit a closed-form solution in terms of matrix operations but numerical
techniques for obtaining the solution are available [see \citeauthor
{MR1997753} (\citeyear{MR1997753}), Chapter 2]. Once the
matrix $v_1$ is known, the linear equation \eqref{vzeq} yields a unique
solution for $v_0$ and $\lambda$ is quadratic in
$v_0$ and linear in $v_1$.

Observe that Theorems \ref{teolin} and \ref{Tpdemodel} characterize the
candidate pair $(\hatp,\mfnct)$, and allow to find the finite-horizon bounds,
but do not address long-run optimality. This stronger property, in fact,
holds only under parameter restrictions and is now studied in detail
for a
single state.
In this case, the linear diffusion yields an extension of the models in
\citet{kim1996dnp} and \citet{wachter2002pac},
%
%
\begin{equation}\label{lindifko}
\cases{
dR_t= (\sigma\nu_0+b\sigma\nu_1Y_t) \,dt+\sigma\,dZ_t,\vspace
*{2pt}\cr
dY_t= -b Y_t \,dt+dW_t,\vspace*{2pt}\cr
d\langle R,Y\rangle_t = \rho\,dt,\vspace*{2pt}\cr
r(Y_t) = r_0.}
\end{equation}
The constants are the same as in \eqref{lindif} except that here $\mu_0
= \sigma\nu_0$, $\mu_1 =
b\sigma\nu_1$ where $\nu_0,\nu_1 \in\reals^n$, for ease of
notation. Note
that $a=1$ and $r_1 = 0$. The Riccati
equation from \eqref{vueq} is
\[
\delta^{-1}v_1^2 + 2b(1+q\rho'\nu_1)v_1-qb^2\nu_1'\nu_1 =
0
\]
for as in (\ref{Edeltadef}). Under Assumption
\ref{assmodlindif}, for $p<0$, the solution $v_1,v_0,\lambda$ from
Theorem \ref{teolin} is
\begin{eqnarray}\label{Eouoptans}
v_1 &= &\delta b\bigl(\sqrt{\Theta}-(1+q\rho'\nu_1)\bigr),\\
v_0 &=&
q\delta\rho'\nu_0 - \frac{1}{\sqrt{\Theta}}\bigl(q\nu_1'\nu_0 +
q\delta
\rho'\nu_0(1+q\rho'\nu_1)\bigr),\\
\lambda&=& pr_0 - \frac{1}{2}q\nu_0'\nu_0 +
\frac{1}{2}\delta^{-1}v_0^2 - qv_0\rho'\nu_0 - \frac{1}{2}v_1,
\end{eqnarray}
where
%
%
\begin{equation}\label{EThetaval}
\Theta=(1+q\rho'\nu_1)^2 + \delta^{-1}q\nu_1'\nu_1.
\end{equation}
The candidate long-run optimal pair $(\hatp,\mfnct)$ is affine in
the state variable
%
%
\begin{eqnarray}\label{Egenouoptimalstrategy}
\pi(y) &=& \frac{1}{1-p}\Sigma^{-1}\bigl(\mu(y) + v_0\sigma\rho- v_1
y\sigma\rho\bigr),
\nonumber
\\[-8pt]
\\[-8pt]
\nonumber
\mfnct(y) &=&v_0 - v_1 y
\end{eqnarray}
and the dynamics of $(Y,R)$ under the candidate long-run martingale
measure are
%
%
\begin{equation}\label{Egenouoptimalelmmdynamics}
\cases{
dR_t = \sigma\,dZ_t,\vspace*{2pt}\cr
dY_t =
\bigl(-b Y_t -\rho'\sigma^{-1}\mu+
(1 - \rho'\rho)(v_0 - v_1Y_t)\bigr)\,dt +dW_t.}
\end{equation}
This pair $(\hatp,\mfnct)$ is indeed long-run optimal but only
under a parameter restriction.
\begin{proposition}\label{Pgeneralou}
Let Assumption \ref{assmodlindif} hold and let $p<0$, $(\hatp,\mfnct
)$ from~\eqref{Egenouoptimalstrategy}
is long-run optimal if
%
%
\begin{equation}\label{Egeneraloucond}
(1-2q\rho'\rho)\sqrt{\Theta}+(1+q\rho'\nu_1)>0.
\end{equation}
%
\end{proposition}

In the case $\nu_1 = -\kappa\rho$ for $\kappa>0$, which still nests the
models of \citet{kim1996dnp} and \citet{wachter2002pac}, the
parameter restriction in \eqref{Egeneraloucond} simplifies as follows.
%
\begin{corollary}\label{Cfadstypeou}
Let Assumption \ref{assmodlindif} hold. For $p < 0$ and $\nu
_1=-\kappa\rho$ for $\kappa\in\reals$. If
$0<q\rho'\rho\leq1/4$ then long-run optimality holds for all
$\kappa$. For $1/4 <q\rho'\rho<1$
long-run optimality holds if
%
%
\begin{equation}\label{Efadstypeoucond}
\kappa< \frac{2}{4q\rho'\rho-1}.
\end{equation}
\end{corollary}

Thus long-run optimality requires a joint restriction on
preferences ($q$) and price dynamics ($\rho'\rho$ and $\kappa$).
First, since $q \rho'\rho<1$, long-run optimality always holds if
$\kappa< \frac{2}{3}$, that is, if risk premia have low covariation
with changes in state variables. If this condition is not
satisfied, long-run optimality still holds regardless of the level
of incompleteness ($\rho'\rho$) if risk aversion is sufficiently low
($q<\frac14$).\vspace*{-1pt} Conversely, if the market is sufficiently incomplete
($\rho'\rho<\frac14$), the restriction holds regardless of
preferences. Hence, a violation of long-run optimality requires a
high sensitivity of risk premia, high risk aversion and a nearly
complete market.

When long-run optimality fails, it does so at different scales,
depending on parameters. The next proposition studies this phenomenon
in the case $\kappa=1$, which corresponds to a continuous time version
of the model of \citet{summers1986dsm}.
\begin{proposition}\label{PfadsisOK}
Let Assumption \ref{assmodlindif} hold. For $p<0$ and $\kappa=1$
from Corollary \ref{Cfadstypeou},
long-run optimality holds if $q\rho'\rho<\frac{3}{4}$. If
$q\rho'\rho\geq\frac{3}{4}$, long-run optimality fails. In
particular:
\begin{longlist}[(iii)]
\item[(i)]
if $q\rho'\rho> \frac{3}{4}$, there exists a finite $T$ such that
$\frac1pE[(X^{\hatp}_T)^p]=-\infty$;

\item[(ii)]
if $q\rho'\rho=\frac{3}{4}$ and $\nu_0=0$, the certainty equivalent
loss is bounded;

\item[(iii)]
if $q\rho'\rho=\frac{3}{4}$ and $\nu_0\ne0$, the certainty
equivalent loss diverges to $\infty$.
\end{longlist}
\end{proposition}

\subsubsection{Calibration}
A calibration to real data shows that long-run optimality holds for typical
levels of risk aversion, in the model with one asset and one state considered
by \citet{barberis2000ilr} and \citet{wachter2002pac}. The state
variable represents the dividend yield and the asset is an equity
index. In
the notation of this section, they use the set of parameter values (in
monthly units) $\rho=-0.935$, $r=0.14\%$, $\sigma=4.36\%$, $\nu_0=0.0788$,
$\kappa=0.8944$, $b=0.0226$. Then, condition \eqref{Efadstypeoucond} is
satisfied for $p>-12.4$, that is, for risk-aversion less than
$13.4$.

%
\begin{figure}

\includegraphics{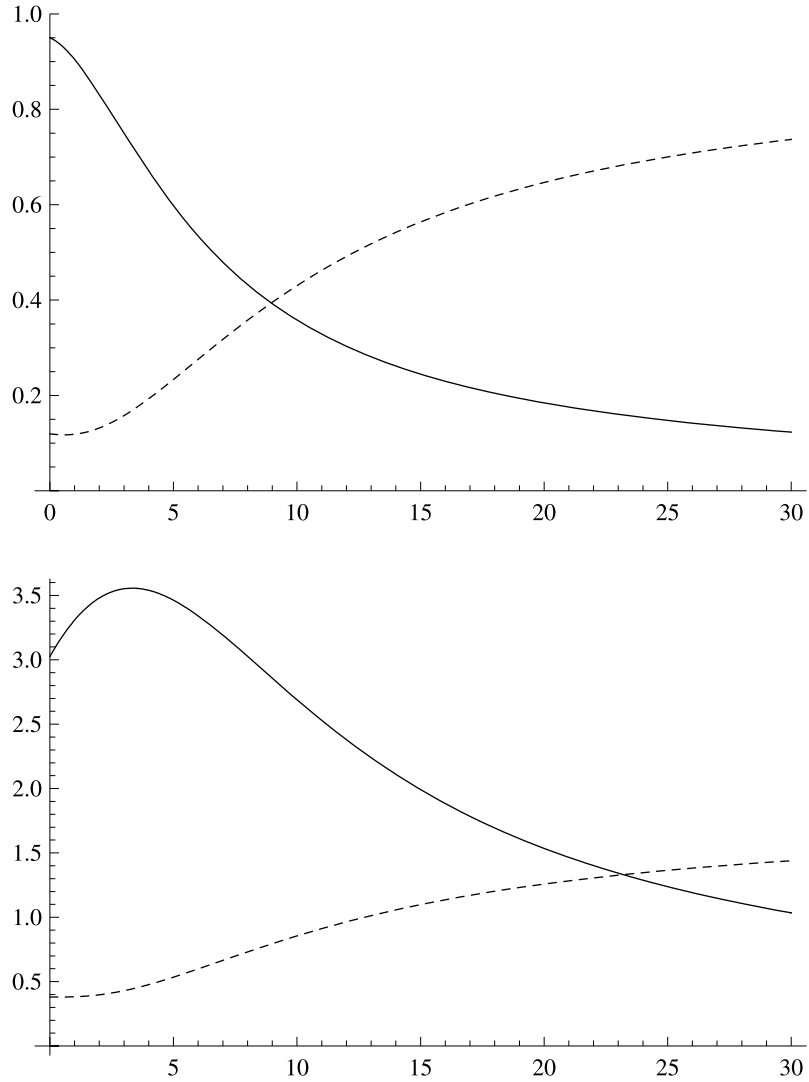}

\caption{Annualized certainty equivalent loss bound (in
percentage points) as a function of the horizon (in years) of the myopic
component (dashed line) and long-run optimal (solid line) portfolios.
Both plots are obtained by equation \protect\eqref{celbou} setting
$\eta$ equal
to the
long-run optimal risk premium, and $\pi$ equal to the long-run optimal
portfolio (solid line) and to the myopic portfolio
$\frac{1}{1-p}\Sigma^{-1}\mu$ (dashed line).
Risk aversion is equal to two ($p=-1$, top) and to five ($p=-4$, bottom).}
\label{longmyopic}
\vspace*{-3pt}
\end{figure}

Figure \ref{longmyopic} compares the finite-horizon performance of the
long-run optimal portfolio to the one of its myopic component. The
plots show the estimates of the correponding upper bounds in \eqref{celbou}:
the myopic component prevails in the short run but its performance
progressively deteriorates as the horizon increases.
The break-even horizon significantly increases with risk-aversion,
passing from nine years for a risk-aversion of two, to twenty-three
years for a risk-aversion of five. Also, the magnitude of the certainty
equivalent loss increases with risk aversion: the differences are
within one percentage point for a risk-aversion of two but increase to
three precentage points for a~risk-aversion of five.

\subsection{Stochastic drifts, volatilities and interest rates}\label{SCIR}
The next model features a single state variable following the
square-root diffusion of \citet{MR0054814}, which simultaneously
affects the interest rate [\citet{MR785475}], the volatilities of
risky assets and their drifts. Note that the model is neither
affine nor quadratic (due to the presence of the term with~$\nu_0$)
and yet the long-run solution admits a simple expression
%
%
\begin{equation}\label{Ecirsetup}
\cases{
dR_t = (\sigma\nu_0 + \sigma\nu_1Y_t)\,dt + \sqrt{Y_t}\sigma
\,dZ_t,\vspace*{2pt}\cr
dY_t = b(\theta- Y_t)\,dt + a\sqrt{Y_t}\,dW_t,\vspace*{2pt}\cr
d\langle R,Y\rangle_t = \rho\,dt,\vspace*{2pt}\cr
r(Y_t) = r_0+r_1 Y_t,}
\end{equation}
where $\sigma\in\reals^{n\times n}$; $\nu_0, \nu_1\in\reals^n$;
$b,\theta, a\in\reals$, $\rho\in\reals^n$ and $r_0,r_1\in\reals
$. The
parameters satisfy the following.
\begin{assum}\label{assmodcir}
$b,\theta,a,r_1\ge0$ and $b\theta> \frac{1}{2}a^2$.
\end{assum}

Assumption \ref{assmodcir} ensures that the state variable $Y$ remains
strictly positive, thereby satisfying Assumption \ref{assmod} with
$E=(0,\infty)$. Guessing a form of the solution $v(y) = v_0\log y + v_1
y$, the
main ODE \eqref{EmainvPDE} becomes
\begin{eqnarray*}
&&pr_0 + pr_1y -
\frac{q}{2}(\nu_0'+y\nu_1')\frac{1}{y}(\nu_0+\nu_1
y) +
\frac{1}{2}\biggl(\frac{v_0}{y}+v_1\biggr)a^2\delta^{-1}y\biggl
(\frac
{v_0}{y}+v_1\biggr)\\
&&\qquad{}+\biggl(\frac{v_0}{y}+v_1\biggr)\bigl(b\theta- qa\rho'\nu_0-(b
+ qa\rho'\nu_1)y\bigr) - \frac{1}{2}a^2\frac{v_0}{y} =
\lambda.
\end{eqnarray*}
Multiplying the above equation by $y$ and setting the constant, linear and
quadratic terms to zero leads to four candidate solutions,
corresponding to any combination of signs in the terms $\pm\sqrt
{\Theta
}$ and $\pm\sqrt{\Lambda}$ below:
\begin{eqnarray*}
v_1 &= &\frac{\delta}{a^2}\bigl(b+qa\rho'\nu_1 \pm
\sqrt{\Theta}\bigr),\\
v_0 &=&
\frac{\delta}{a^2}\bigl(-(c-qa\rho'\nu_0)
\pm\sqrt{\Lambda}\bigr),\\
\lambda&=& pr_0 - q\nu_0'\nu_1 + \frac{a^2}{\delta
}v_0v_1 -
v_0(b+qa\rho'\nu_1) +
v_1(b\theta-qa\rho'\nu_0)
\end{eqnarray*}
with
%
%
\begin{eqnarray}\label{EcirThetaLambdaval}
c&=&b\theta- \tfrac{1}{2}a^2,\nonumber\\
\Theta&=& (b+qa\rho'\nu_1)^2 +
\frac{a^2}{\delta}(q\nu_1'\nu_1 - 2pr_1),\\
\Lambda&=& (c -qa\rho'\nu_0)^2 +
\frac{a^2}{\delta}q\nu_0'\nu_0.\nonumber
\end{eqnarray}
When $p < 0, r_1 > 0$,
%
%
\begin{eqnarray}\label{Ecompare}
\Theta&>& (b+qa\rho'\nu_1)^2>0,
\nonumber
\\[-8pt]
\\[-8pt]
\nonumber
\Lambda&>& (c-qa\rho'\nu_0)^2>0.
\end{eqnarray}
Under $\hat{P}$, $Y$ has again four possible dynamics:
\[
dY_t = \bigl(\tfrac{1}{2}a^2\pm\sqrt{\Lambda}
\pm\sqrt{\Theta}Y_t\bigr)\,dt + a\sqrt{Y_t}\,dW_t.
\]
The choice of $-\sqrt{\Theta}$ and $\sqrt{\Lambda}$ ensures that $Y$
satisfies Assumption \ref{assmod} under $\hat{P}^y$ and is $\hat{P}^y$-tight
in $(0,\infty)$ for each $y\in(0,\infty)$. The latter statement
follows by
the positivity of $\Theta,\Lambda$ [see \citet{MR1326606}, Corollary
5.1.11, and
the discussion immediately after Assumption \ref{boundedpotential}].
Thus the candidate optimizer is
%
%
\begin{eqnarray}\label{Eciroptans}
v_1 &=& \frac{\delta}{a^2}\bigl(b+qa\rho'\nu_1-
\sqrt{\Theta}\bigr),\nonumber\\
v_0 &=&
\frac{\delta}{a^2}\bigl(-(c-qa\rho'\nu_0)
+ \sqrt{\Lambda}\bigr),\\
\lambda&=& pr_0 - q\nu_0'\nu_1 + \frac{a^2}{\delta
}v_0v_1 -
v_0(b+qa\rho'\nu_1) +
v_1(b\theta-qa\rho'\nu_0).\nonumber
\end{eqnarray}
The candidate long-run optimal policies $\lrpair$ are
%
\[
\hatp(y)= \frac{1}{1-p}{\Sigma}^{-1}\bigl(\mu(y) + \sigma\rho
a(v_0 + v_1 y)\bigr) ,\qquad\eta(y)=
\frac{v_0}{y}+v_1
\]
and the candidate long-run martingale measure is
\[
\cases{
dR_t =
\sqrt{Y_t}\sigma\,dZ_t,\vspace*{2pt}\cr
dY_t = \bigl(
b(\theta-Y_t) - qa(\rho'\nu_0+\rho'\nu_1 Y_t)+
a^2(1-\rho'\rho)(v_0 + v_1 Y_t)\bigr)\,dt \vspace*{2pt}\cr
\hspace*{30pt}{}+ a\sqrt{Y_t}\,dW_t.}
\]
Long-run optimality obtains under the following conditions.
\begin{proposition}\label{LcirisOK}
Let Assumption \ref{assmodcir} hold. For $p<0$, long-run optimality
holds if
%
%
\begin{eqnarray}\label{parcir}
(1-2q\rho'\rho)\sqrt{\Lambda}
+(c-qa\rho'\nu_0) &> &0,
\nonumber
\\[-8pt]
\\[-8pt]
\nonumber
(1-2q\rho'\rho)\sqrt{\Theta} +
(b+qa\rho'\nu_1) &>& 0.
\end{eqnarray}
\end{proposition}

The main economic message of this parametric restriction is the same as
in the
previous example. Long-run optimality holds if either one of the following
conditions is satisfied: the covariation of risk premia with state
shocks is
small ($a \rho'\nu_0, a\rho'\nu_1 \approx0$), the market is sufficiently
incomplete ($\rho'\rho\ll1$) or risk aversion is low ($1-p\ll\infty$).

\section{Conclusion}\label{secconc}
Long-run analysis is a tractable and yet nontrivial framework for
dynamic portfolio choice and derivatives pricing in incomplete markets,
and yields simple expressions for portfolios and risk premia.
Long-run policies admit closed-form solutions even in cases in which
finite-horizon policies do not, and the finite-horizon performance of
long-run policies has a simple expression.

Long-run optimality entails that the certainty equivalent loss vanishes
for long horizons and requires some joint restrictions on preferences
and asset dynamics. It does not hold at the intersection of three
extreme cases: risk premia highly
co-varying with state shocks, a nearly complete market and high risk
aversion. Otherwise, long-run optimality holds and time-homogeneous
portfolios are approximately optimal for long horizons.

\begin{appendix}\label{app}

\section{\texorpdfstring{Proof of Section \lowercase{\protect\ref{secmodel}}}{Proof of Section 1}}


\begin{pf*}{Proof of Lemma \protect\ref{LHolderdualitybound}}
Denote by $q=\frac{p}{p-1}$.\vspace*{-1pt} In the case $0<p<1$, H\"{o}lder's
inequality with $\tilde{p} = \frac{1}{p}$ and $\tilde{q} =
\frac{\tilde{p}}{\tilde{p}-1} = \frac{1}{1-p}$ yields
\begin{eqnarray*}
\espalt{P}{X^p} &= &\espalt{P}{(XM)^pM^{-p}}\leq
\espalt{P}{(XM)^{p\tilde{p}}}^{1/\tilde{p}}\espalt{P}{M^{-p\tilde
{q}}}^{1/\tilde{q}}\\
&= &\espalt{P}{XM}^{1/\tilde{p}}\espalt{P}{M^q}^{1-p}\leq
\espp{M^q}^{1-p},
\end{eqnarray*}
because $\espalt{P}{XM}\leq1$, and the claim follows dividing by
$p>0$. If $p < 0$ ($0 <  q < 1$) H\"{o}lder's inequality with $\tilde{p}
= \frac{1}{1-q}$, $\tilde{q} = \frac{1}{q}$ yields
\begin{eqnarray*}
\espalt{P}{M^q}^{1-p} &=& \espalt{P}{(XM)^qX^{-q}}^{1-p}\leq
\espalt{P}{(XM)^{q\tilde{q}}}^{{(1-p)}/{\tilde{q}}}\espalt
{P}{X^{-q\tilde{p}}}^{{(1-p)}/{\tilde{p}}}\\
&=& \espalt{P}{XM}^{-p}\espalt{P}{X^p}\leq\espp{X^p}
\end{eqnarray*}
and the claim now follows dividing by $p<0$. In both cases, the
inequality becomes an equality when $\espalt{P}{X M}=1$ and
$X^{p-1}$ is proportional to $M$.
\end{pf*}


\section{\texorpdfstring{Proofs of Section \lowercase{\protect\ref{seclongrun}}}{Proofs of Section 2}}

\begin{pf*}{Proof of Theorem \protect\ref{Tpdemodel}}
Since the Brownian motions $Z$ and $W$ are partially correlated, the
following orthogonal decomposition holds (see the discussion following
Assumption \ref{assmod}):
%
%
\begin{equation}\label{ortdec}
dZ_t=\rho(Y_t) \,dW_t+\bar\rho(Y_t) \,dB_t,
\end{equation}
where $B=(B^1,\ldots,B^n)$ is a $n$-dimensional Brownian motion
independent of $W$ and the matrix $\bar\rho(y)$ is defined by the
identity $(\rho\rho')(y)+(\bar\rho\bar\rho')(y)=I_n$.
For $t\ge0$, define the process $D$ by
%
%
\begin{eqnarray}\label{Ephatprnderiv}
D_t& =& \mathcal{E}\biggl(\int_0^{\cdot}\bigl(-q\ups'\Sigma
^{-1}\mu+
(A-q\Upsilon'\Sigma^{-1}\Upsilon)\nabla v\bigr)'(a')^{-1}\,dW_t
\nonumber
\\[-8pt]
\\[-8pt]
\nonumber
&&\hspace*{72pt}\quad{}-
\int_0^{\cdot}q(\Sigma^{-1}\mu+ \Sigma^{-1}\ups\nabla v
)'\sigma\bar{\rho} \,dB_t\biggr)_t.
\end{eqnarray}
By assumption, the operators associated to the models $P$ and $\hat{P}$
satisfy Assumption~\ref{assmod}. Since Assumption~\ref{regmod} also
holds, \citeauthor{MR2152242} [(\citeyear{MR2152242}), Theorem 2.4, Remark 2.4.2] implies that for
each $y\in E$ and $t\geq0$,~$P^y$ and
$\hat{P}^y$ are equivalent on $\mathcal{F}_t$ with
%
%
\begin{equation}\label{Ernderivffilt}
\frac{d\hat{P}^y}{dP^y}\bigg|_{\mathcal{F}_t} = D_t.
\end{equation}
Thus, $(D_t,\mathcal{F}_t)_{t\geq0}$ is a $P^y$-martingale. With this
notation, it suffices to prove that given the solution pair $(v,\lambda
)$ to \eqref{EmainvPDE}, with $\pi,\eta$ as in \eqref{lrport}, the
following~$P^y$ almost sure identities hold:
\begin{eqnarray}
\label{priestas}
(X^{\hatp}_T)^p &=& e^{\lambda T+v(y)-v(Y_T)}D_T,\\
\label{duaestas}
(M^{\mfnct}_T)^q &=& e^{({1}/{(1-p)})(\lambda
T+v(y)-v(Y_T))}D_T.
\end{eqnarray}
Indeed, if \eqref{priestas} and \eqref{duaestas} hold, then \eqref{priest}
and \eqref{duaest} follow by taking expectations with respect to $P^y$.
Consider first \eqref{priestas}.
Passing to logarithms, it suffices to prove that
%
%
\begin{equation}\label{Epritemp1}
p\log X^{\pi}_T - \log D_T = \lambda T + v(y) - v(Y_T).
\end{equation}
The first term on the left-hand side of \eqref{Epritemp1} is
\[
p\log X^{\pi}_T = \int_0^T\biggl(pr + p\pi'\mu-
\frac{1}{2}p\pi'\Sigma\pi\biggr)\,dt +
\int_0^{T} p\pi'\sigma\,dZ_t.
\]
Substituting $\hatp= \frac{1}{1-p}\Sigma^{-1}(\mu+ \ups\nabla
v)$, the decomposition $Z=\rho W + \bar{\rho} B$ and collecting terms
\begin{eqnarray}\label{Epritemp2}
p\log X^{\pi}_T &=& \int_0^T\biggl(pr - \frac
{1}{2}q(1+q)\mu
'\Sigma
^{-1}\mu-
q^2\mu'\Sigma^{-1}\Upsilon\nabla v \nonumber\\
&&\hspace*{58pt}\qquad{}+ \frac{1}{2}q(1-q)\nabla
v'\Upsilon'\Sigma^{-1}\Upsilon\nabla v\biggr)\,dt\\
&&{}-q\int_0^T(\mu+ \Upsilon\nabla v)'\Sigma^{-1}\sigma\rho
\,dW_t-q\int_0^T(\mu+ \Upsilon\nabla v)'\Sigma^{-1}\sigma\bar
{\rho} \,dB_t.\nonumber
\end{eqnarray}
The second term in the left-hand side of \eqref{Epritemp1} follows
from \eqref{Ephatprnderiv},
\begin{eqnarray}\label{Epritemp4}
\log D_T &= &\int_0^T\biggl(-\frac{1}{2}q^2\mu'\Sigma^{-1}\mu+
q(1-q)\mu'\Sigma^{-1}\Upsilon\nabla v \nonumber\\[-2pt]
&&\hspace*{8pt}\qquad{}- \frac{1}{2}\nabla
v'\bigl(A+(q^2-2q)\Upsilon'\Sigma^{-1}\Upsilon\bigr)\nabla v
\biggr)\,dt
\nonumber
\\[-10pt]
\\[-10pt]
\nonumber
&&{}+\int_0^T\bigl(\nabla v'a - q(\mu+ \Upsilon\nabla
v)'\Sigma^{-1}\sigma\rho\bigr)\,dW_t\\[-2pt]
&&{} - q\int_0^T(\mu+
\Upsilon\nabla v)\Sigma^{-1}\sigma\bar{\rho}\,dB_t.\nonumber
\end{eqnarray}
Subtracting \eqref{Epritemp4} from \eqref{Epritemp2} yields
\begin{eqnarray*}
&&p\log X^{\pi}_T-\log D_T\\[-2pt]
&&\qquad=
\int_0^T\biggl(pr - \frac{1}{2}q\mu'\Sigma^{-1}\mu
- q\mu'\Sigma^{-1}\Upsilon\nabla v + \frac{1}{2}\nabla
v'(A-q\Upsilon'\Sigma^{-1}\Upsilon)\nabla v\biggr)\,dt
\\[-2pt]
&&\qquad\quad{}-\int_0^T \nabla v'a \,dW_t.
\end{eqnarray*}
Now, It\^{o}'s formula allows us to substitute
%
%
\begin{equation}\label{Eitov}
\quad -\int_0^T\nabla v'a\,dW_t = v(y) - v(Y_T) + \int_0^T\nabla v'b \,dt +
\frac{1}{2}\int_0^T \trace(AD^2v)\,dt
\end{equation}
and the claim \eqref{Epritemp1} follows by recalling that $(v,\lambda
)$ solves \eqref{EmainvPDE}.

Consider now the equality \eqref{duaestas}. Again, by taking
logarithms it suffices to show that
%
%
\begin{equation}\label{Eduatemp1}
q\log M^{\eta}_T - \log D_T = \frac{1}{1-p}\bigl(\lambda T +
v(y)-v(Y_T)\bigr).
\end{equation}
The first term in the left-hand side is equal to (plugging $\mfnct=
\nabla v$ and $Z=\rho W + \bar{\rho} B$)
\begin{eqnarray*}
q\log M^{\mfnct}_T &=& -q\int_0^T r \,dt + q\log\mathcal{E}\biggl
(-\int
_0^{\cdot}(\Sigma^{-1}\mu+ \Sigma^{-1}\ups\mfnct)'\sigma
\,dZ_t\\[-2pt]
&&\hspace*{174pt}{} + \int_0^{\cdot}\mfnct'a\,dW_t\biggr)_T\\[-2pt]
&=& -q\int_0^T\biggl(r +\frac{1}{2}\mu'\Sigma^{-1}\mu+\frac
{1}{2}\nabla
v'(A-\Upsilon'\Sigma^{-1}\Upsilon)\nabla v\biggr)\,dt\\[-2pt]
&&{}+ q\int_0^T\bigl(\nabla v'a -(\mu+ \Upsilon\nabla
v)'\Sigma^{-1}\sigma\rho\bigr) \,dW_t\\[-2pt]
&&{} - q\int_0^T(\mu+
\Upsilon\nabla v)\Sigma^{-1}\sigma\bar{\rho}\,dB_t
\end{eqnarray*}
and subtracting \eqref{Epritemp4} yields
\begin{eqnarray*}
q\log M^{\mfnct}_T-\log D_T &=&
\frac{1}{1-p}\int_0^T\biggl(pr - \frac{1}{2}q\mu
'\Sigma
^{-1}\mu
- q\mu'\Sigma^{-1}\Upsilon\nabla v \\
&&\hspace*{76pt}{}+ \frac{1}{2}\nabla
v'(A-q\Upsilon'\Sigma^{-1}\Upsilon)\nabla v\biggr)\,dt\\
&&{}-\frac1{1-p}\int_0^T \nabla v'a \,dW_t.
\end{eqnarray*}
As in the previous case, \eqref{Eduatemp1} now follows by
substituting \eqref{Eitov} and recalling that $(v,\lambda)$ solves
\eqref{EmainvPDE}.\qed

\noqed\end{pf*}
%

The proof of Theorem \ref{Tlongrunopt} requires two lemmas.
\begin{lemma}\label{Ltightnesslemma}
Let $\phi\in C(E;\reals)$, $\phi> 0$ and let $(\mu_T)_{T\geq0}$ be
a tight
family of probability measures on $(E,\borel(E))$. Then
%
%
\begin{equation}\label{Etightnessresult}
\liminf_{T\uparrow\infty} \frac{1}{T}\log\int\phi\,d\mu_T \geq0.
\end{equation}
\end{lemma}
\begin{pf}
By monotonicity, it suffices to prove the result for $\phi$ bounded. Let
$(t_n)_{n\ge1}$ be an increasing sequence satisfying
$t_n\uparrow\infty$ and
%
%
\begin{equation}\label{Etempeqn}
\liminf_{T\uparrow\infty}\frac{1}{T}\log\int_E\phi\,d\mu_T
=\lim_{n\conv\infty}\frac{1}{t_n}\log\int_E \phi\,d\mu_{t_n}.
\end{equation}
Since $\phi$ is bounded above, the limit in \eqref{Etempeqn} is
nonpositive. Since the measures $(\mu_T)_{T\geq0}$ are
tight, they
are relatively compact with respect to the topology of weak
convergence. Thus, up to a subsequence, there exists a probability
measure~$\mu$ on~$E$ such that
\[
\lim_{n\conv\infty} \int_E \phi\,d\mu_{t_n} = \int_E \phi\,d\mu
\in
(0,\infty)
\]
because $\phi$ is continuous, bounded and positive. Thus, for any
$0 < \eps< \int_E\phi\,d\mu$, there is a $n_{\eps}$ such that
$n\geq
n_{\eps}$ implies
\[
\frac{1}{t_n}\log\int_{E}\phi\,d\mu_{t_n} \geq\frac{1}{t_n}\log
\biggl(\int_E\phi
\,d\mu- \eps/2\biggr).
\]
Hence, taking $n\uparrow\infty$ it follows that the limit in
\eqref{Etempeqn} is indeed $0$.
\end{pf}


\begin{lemma}\label{Lgeneratorupperbound}
Let $(Q^y)_{y\in E}$ be a solution to the martingale problem for
the operator $L$ on $E$, where
\[
L=\frac{1}{2}\sum_{i,j=1}^{k}\Sigma(y)^{ij}\frac{\partial
^2}{\partial
y_i\,\partial
y_j} + \sum_{i=1}^{k}\mu(y)^{i}\frac{\partial}{\partial y_i}.
\]
If $f\in C^2(E,\reals)$, $f>0$, then $E_{Q}^y[f(Y_T)]\leq f(y) +
(0\vee\sup_E Lf)T$.
\end{lemma}

\begin{pf}
Let $f\in C^2(E,\reals)$, $f>0$. Since $(Q^y)_{y\in E}$ solves the martingale
problem for $L$ on $E$, the process
\[
f(Y_T) - \int_0^T(Lf)(Y_t)\,dt
\]
is a local martingale under $Q^y$. Let $(\tau_n)_{n\ge1}$ be a
reducing sequence of stopping times for this local martingale. Then
\[
E_{Q}^y[f(Y_{T\wedge\tau_n})]=f(y) +
E_{Q}^y\biggl[\int_0^{T\wedge\tau_n}(Lf)(Y_t)\,dt\biggr]\leq f(y)
+ \Bigl(0\vee
\sup_{E}(L f)\Bigr)T
\]
and the thesis follows by Fatou's lemma, taking $n\uparrow\infty$.\vspace*{-2pt}
\end{pf}


\begin{pf*}{Proof of Theorem \protect\ref{Tlongrunopt}}
Applying the inequality in \eqref{duaineq} to equations~\eqref
{priest} and~\eqref{duaest} from Theorem \ref{Tpdemodel} gives
%
%
\begin{eqnarray}\label{Eliminfsupvalues}
\qquad0&\leq&\liminf_{T\conv\infty} \frac1p \biggl(
\frac1T\log E_P^y[(M^{\mfnct}_T)^q]^{1-p}-\frac1T\log
E_P^y[(X^{\hatp}_T)^p]
\biggr)\nonumber\\
\qquad&\leq&\limsup_{T\conv\infty}
\frac1{pT}\log E_P^y[(M^{\mfnct}_T)^q]^{1-p}-
\liminf_{T\conv\infty}
\frac1{pT}\log E_P^y[(X^{\hatp}_T)^p]\\
\qquad&=&\limsup_{T\conv\infty}\frac{1-p}{pT}\log E_{\hat
{P}}^y\bigl[e^{- ({1}/{(1-p)})v(Y_T)}\bigr]
- \liminf_{T\conv\infty}\frac{1}{pT}\log
E_{\hat{P}}^y\bigl[e^{-v(Y_T)}\bigr].\nonumber
\end{eqnarray}
Thus it is sufficient to prove for $0 < p < 1$ that
\begin{eqnarray}
\label{lopos}
\limsup_{T\conv\infty}\frac{1}{T}\log E_{\hat P}^y\biggl[\exp
\biggl(-\frac
{1}{1-p}v(Y_T)\biggr)\biggr]&\leq&0,\\
\label{uppos}
\liminf_{T\conv\infty}\frac{1}{T}\log E_{\hat P}^y[\exp
(-v(Y_T))]&\geq&0
\end{eqnarray}
and for $p < 0$ that
\begin{eqnarray}
\label{upneg}
\liminf_{T\conv\infty}\frac{1}{T}\log E_{\hat P}^y\biggl[\exp
\biggl(-\frac
{1}{1-p}v(Y_T)\biggr)\biggr]&\geq&0,\\
\label{loneg}
\limsup_{T\conv\infty}\frac{1}{T}\log E_{\hat P}^y[\exp
(-v(Y_T))]&\leq&0
\end{eqnarray}
for $p<0$. The lower bounds \eqref{upneg} and \eqref{uppos} follow
from the
application of Lemma~\ref{Ltightnesslemma} to the functions $\phi
=\exp(-\frac{1}{1-p}v)$ and $\phi= \exp(-v)$, respectively.
For the upper bounds, first denote by
\[
Lf = \nabla f'\bigl(b - q\ups'\Sigma^{-1}\mu+
(A-q\Upsilon'\Sigma^{-1}\Upsilon)\nabla v \bigr) +
\tfrac{1}{2}\trace(AD^2f)
\]
and observe that, for any $\alpha\in\reals$,
\begin{eqnarray*}
L(e^{\alpha v}) &= &\alpha e^{\alpha v}\biggl(\nabla v'
\bigl(b-q\ups'\Sigma^{-1}\mu+ (A-q\Upsilon'\Sigma^{-1}\Upsilon
)\nabla v\bigr) \\
&&\hspace*{89pt}{}+ \frac{1}{2}\trace(AD^2v) + \frac
{1}{2}\alpha\nabla v'A\nabla v\biggr)\\
&=& \alpha e^{\alpha v}\biggl(\frac{1}{2}\nabla v'\bigl((1+\alpha)A -
q\Upsilon'\Sigma^{-1}\Upsilon\bigr)\nabla v + \lambda- {pr} + \frac
q2\mu'\Sigma^{-1}\mu\biggr),
\end{eqnarray*}
where the second equality follows from \eqref{EmainvPDE}.
For $p < 0$, consider $\alpha= -1$,
\[
L(e^{-v}) = e^{-v}\biggl(\frac q2\nabla v'\Upsilon'\Sigma
^{-1}\Upsilon\nabla v - \lambda+ pr - \frac q2\mu
'\Sigma^{-1}\mu\biggr).
\]
Assumption (ii) implies that the right-hand side is bounded by some
constant~$K$ and Lemma \ref{Lgeneratorupperbound} yields
\[
E_{\hat P}^y\bigl[e^{-v(Y_T)}\bigr]\leq e^{-v(y)} + (K\vee0)T
\]
and hence,
\[
\limsup_{T\conv\infty}\frac{1}{T}\log E_{\hat P}^y\bigl
[e^{-v(Y_T)}\bigr]\leq0.
\]
Similarly, for $0 < p < 1$ consider $\alpha= -\frac{1}{1-p}$,
\begin{eqnarray*}
L\bigl(e^{-({1}/{(1-p)})v}\bigr) &=& \frac{1}{1-p}e^{-
({1}/{(1-p)})v}\\
&&\times{}\biggl(-\frac q2\nabla v'(A-\Upsilon'\Sigma
^{-1}\Upsilon
)\nabla v - \lambda+ pr - \frac q2\mu'\Sigma^{-1}\mu
\biggr).
\end{eqnarray*}
Again, the right-hand side is bounded by $K$ and Lemma
\ref{Lgeneratorupperbound} yields
\[
E_{\hat P}^y\bigl[e^{-({1}/{(1-p)})v(Y_T)}\bigr]\leq e^{-({1}/{(1-p)})v(y)} +
(K\vee0)T
\]
and the claim follows as in the previous case.
\end{pf*}

The proofs of Theorems \ref{TLambdastructure} and \ref{Tbellmansolution}
are obtained by adapting the arguments in \citet{MR2206349} to the present
setting. Because the structure of these proofs remains the same, the arguments
are not repeated in detail, focusing instead on the necessary modifications.

Henceforth, all references in italics point to \citet{MR2206349}. For
ease of notation, set
%
%
\begin{equation}\label{EKSnotation}
\qquad\hat{A} = A-q\Upsilon'\Sigma^{-1}\Upsilon,\qquad\tilde{b} =
b-q\Upsilon'\Sigma^{-1}\mu,\qquad V = pr - \frac
{q}{2}\mu
'\Sigma
^{-1}\mu
\end{equation}
and define the quasilinear operator $M$ acting on $f\in C^2(E)$ by
%
%
\begin{equation}\label{Equasilinop}
Mf = \tfrac{1}{2}\trace(AD^2f) + \tfrac{1}{2}\nabla f'\hat{A}\nabla f +
\tilde{b}'\nabla f + V
\end{equation}
so that \eqref{EmainvPDE} becomes $Mv = \lambda$. The following results
carry over immediately from \citet{MR2206349} and \citet{MR0244627} with
only the obvious notational
changes.
\begin{lemma}[({Lemma $2.4$} and the discussion following)]\label
{LKSThm26Lemmas}
Let Assumption~\ref{regmod} hold. Let $\lambda\in\reals$. Then:
\begin{longlist}[(ii)]
\item[(i)] If, for each $N$ there exists a $v_N\in C^{2,\gamma
}(E_N,\reals)$
satisfying $M v_N = \lambda$ in $E_N$, then there exists a $v\in
C^{2,\gamma}(E,\reals)$ satisfying $Mv =
\lambda$.
\item[(ii)] If $\lambda_{m}\,{\conv}\,\lambda$ and for each $m$ there is a
$v_m\,{\in}\,C^{2,\gamma}(E,\reals)$ solving \mbox{$Mv_m\,{=}\,\lambda_m$}, then there is a
$v\in
C^{2,\gamma}(E,\reals)$ solving $Mv = \lambda$.
\end{longlist}
\end{lemma}

\begin{remark} Note that Lemma $2.4$ requires the uniform
ellipticity of
diffusion matrices in the form of conditions (A1),
(A2). Under
Assumption~\ref{regmod}, these conditions are satisfied on each $E_N$ by
continuity and
positivity. In the present setting, the conclusion of Lemma
2.4 is
that, if $v_N\in C^{2,\gamma}(E_N,\reals)$ solves $Mv_N = \lambda$
in $E_N$,
then for each $n$ there exist constants $A_n,B_n,C_n > 0$ such that for
$N>2n$
\[
\sup_{E_n}|\nabla v_N| \leq\max\{A_n + B_n\lambda,
C_n\}.
\]
\end{remark}

\begin{lemma}[{[Theorem 8.4 in \citet{MR0244627}]}]\label{LLUpde}
Let Assumptions \ref{regmod} and \ref{regmoddomain} hold. For any
$N$, if
there exist functions $f_1,f_2 \in  C^{2,\gamma}(E_{N+1},\reals)\cap
C(\bar{E}_{N+1},\reals)$ satisfying $Mf_1>\lambda$ and $Mf_2<\lambda
$ in
$E_N$, then there exists a function $v_N\in C^{2,\gamma}(E_N,\reals)$
satisfying $Mv_N = \lambda$.
\end{lemma}

\begin{pf*}{Proof of Theorem \protect\ref{TLambdastructure}}
In \citet{MR2206349} Theorem \ref{TLambdastructure} is split into two
theorems: {Theorem $2.6$} shows that $\Lambda= [\lambda
_c,\infty
)$, while
{Theorem $3.2$} shows that $Y$ is $(\hat{P}^y)_{y\in E}$-transient
for $\lambda> \lambda_c$. The following arguments show that the
conclusions of both {Theorem $2.6$} and {Theorem $3.2$}
remain valid.

Regarding {Theorem $2.6$}, given Lemmas \ref{LKSThm26Lemmas} and
\ref{LLUpde}, it suffices to prove that:
\begin{longlist}[(A)]
\item[(A)] There exists a $\lambda_0\in\reals$ such that for each
$N$ there
is a
$v_N\in C^{2,\gamma}(E_N,\reals)$ satisfying $Mv_N = \lambda_0$ in $E_N$.
\item[(B)] The set $\{\lambda\in\reals: \exists v\in C^{2,\gamma
}(E,\reals)
\mbox{ satisfying } Mv \geq\lambda\}$ is unbounded from above.
\item[(C)]$\lambda_c \equiv\inf\Lambda> -\infty$.
\end{longlist}

Indeed, if (A) holds true then Lemma \ref{LKSThm26Lemmas} yields a solution
$v\in C^{2,\gamma}(E,\reals)$ to $Mv=\lambda_0$ proving
$\Lambda\neq\varnothing$. Now, let $\lambda\in\Lambda$, let
$v\in C^{2,\gamma}(E,\reals)$ satisfy \mbox{$Mv = \lambda$} and let $\tilde
{\lambda}
> \lambda$. If (B) holds, then there is a $f_1\in
C^{2,\gamma}(E,\reals)$ satisfying \mbox{$Mf_1>\tilde{\lambda}$}. Thus, with
$f_2=v$, Lemma \ref{LLUpde} applies for any $N$ and hence, by
Lemma~\ref{LKSThm26Lemmas}, it follows that $\tilde{\lambda}\in
\Lambda$.
Thus, $\lambda\in\Lambda\Rightarrow[\lambda,\infty)\subseteq
\Lambda$. By
(C) it follows that $\Lambda$ is bounded from below and by Lemma
\ref{LKSThm26Lemmas} it follows that $\Lambda$ is closed. Thus,
$\Lambda=
[\lambda_c,\infty)$ is the desired result.

Therefore, it remains to prove (A), (B) and (C). First, note that with
$a$, $\sigma$ the unique symmetric, positive definite square roots of
$A$ and
$\Sigma,$ respectively, it follows that
\[
\hat{A}= a(1_{k} -
q\rho'\rho)a = (1-q)aa + qa(1_{k} - \rho'\rho)a.\vadjust{\goodbreak}
\]
By construction $1_n -\rho\rho'\geq0$ and hence, $1_k-\rho'\rho
\geq
0$. Thus, setting
\[
c =
\cases{
1-q, &\quad $p < 0,$\vspace*{2pt}\cr
1, &\quad $0 < p < 10$,
}\qquad
\bar{c} =
\cases{
1, &\quad $p < 0,$\vspace*{2pt}\cr
1-q, &\quad $0 < p < 1,$}
\]
it follows that $c,\bar{c} > 0$ are such that on $E$
%
%
\begin{equation}\label{EAhatArel}
c A\leq\hat{A} \leq\bar{c} A.
\end{equation}
Define the linear
operators $L^{c}$ and $L^{\bar{c}}$ acting on $f\in C^2(E)$ by
\[
L^{c}f = \tfrac{1}{2}\trace(AD^2f) + \tilde{b}'\nabla f + cVf,\qquad
L^{\bar{c}}f = \tfrac{1}{2}\trace(AD^2f) + \tilde{b}'\nabla f + \bar{c}Vf.
\]
Let $\lambda^{*}$ and $\bar{\lambda}^{*}$ denote the generalized
principal eigenvalue for $L^{c}$
and $L^{\bar{c}}$ on~$E$ [\citet{MR1326606}, Chapter 4.3]. Assumptions
\ref{regmod}, \ref{regmoddomain} and \ref{boundedpotential} imply that~%
$\lambda^{*},\allowbreak\bar{\lambda}^{*}\in\reals$. By~\eqref{EAhatArel},
it follows
that
%
%
\begin{equation}\label{EMLbounds}
\frac{1}{cg}L^{c}g\leq Mf\leq\frac{1}{\bar{c}\bar{g}}
L^{\bar{c}}\bar{g},
\end{equation}
where $g = e^{c f}$ and $\bar{g} = e^{\bar{c}f}$. For any $\beta>
\max[\lambda^{*},\bar{\lambda}^{*},1]$, let $\lambda=
\frac{\beta}{2}(\frac{1}{c}+\frac{1}{\bar{c}})$. By
construction there
exist $g,\bar{g}\in C^{2,\gamma}(E)$, $g,\bar{g} > 0$ satisfying
$L^{c}g =
\beta g$ and $L^{\bar{c}}\bar{g} = \beta\bar{g}$ in $E$. Set $f_1 =
\frac{1}{c}\log g$. By \eqref{EMLbounds}, it follows that on $E$
\[
Mf_1-\lambda\geq\frac{1}{cg}L^cg -\lambda=
\frac{\beta}{2}\biggl(\frac{1}{c}-\frac{1}{\bar{c}}\biggr) > 0.
\]
Similarly, setting $f_2 = \frac{1}{\bar{c}}\log\bar{g}$ and using
\eqref{EMLbounds} it follows that on $E$
\[
Mf_2 - \lambda\leq\frac{1}{\bar{c}\bar{g}}L^{\bar{c}}\bar{g} -
\lambda=
\frac{\beta}{2}\biggl(\frac{1}{\bar{c}} - \frac{1}{c}\biggr) < 0.
\]
Therefore, (A) holds by first applying Lemma \ref{LLUpde} and then Lemma
\ref{LKSThm26Lemmas}(i) for $\lambda_0 = \lambda$.
Furthermore, since $\beta$ can be taken arbitrarily large, (B) holds as
well. Regarding (C), let $\lambda\in\Lambda$ and let $v$ be the
associated function solving $Mv =
\lambda$ and let $g = e^{c v}$. By \eqref{EMLbounds}
\[
\lambda= \sup_{E}\ (Mv)\geq\frac{1}{c}\sup_{E}\biggl(\frac
{L^{c}g}{g}\biggr)\geq\frac{1}{c}\mathop{\inf_{g\in C^2(E,\reals)}}_{g
> 0}\sup_{E}\biggl(\frac{L^{c}g}{g}\biggr)=\frac{\lambda^{*}}{c},
\]
where the last equality follows by \citet{MR1326606}, Theorem 4.4.5. Thus
$\lambda_c>-\infty$.

Regarding Theorem $3.2$, note that by \eqref{EAhatArel},
Remark
$2.1$ holds for $c$, $\bar{c}$ and the argument in
Lemma $3.1$ carries over exactly, up to obvious changes in
notation.
\end{pf*}

Before proving Theorem \ref{Tbellmansolution}, the following
definitions and
results are needed from the theory of large deviations for occupancy
times of
diffusions. Let Assumptions \ref{regmod}, \ref{regmoddomain} and
\ref{boundedpotential} hold. To make the dependence upon $\lambda$ specific,
for $\lambda\in\Lambda$ let
$(\hat{P}^{\lambda,y})_{y\in E}$ be the measure in Theorem \ref{Tpdemodel}
and let $\hat{L}^{\lambda}$ be the operator associated to
$(\hat{P}^{\lambda,y})_{y\in E}$. Let $M_1(E)$ denote the space
of Borel probability measures on $E$. Define the function
$I^{\lambda}: M_1(E)\mapsto\reals$ by
\[
I^{\lambda}(\mu) = -\inf_{u\in\mathcal{U}}\int_{E}\frac{\hat
{L}^{\lambda
}u}{u}\,d\mu,
\]
where
\[
\mathcal{U} = \biggl\{u\in C^2(E,\reals)\ \Big|\ u(x)\geq\eps_u > 0,
\displaystyle\frac{\hat{L}^{\lambda}u}{u}\mbox{ bounded}\biggr\}.
\]

It is clear that
$I^{\lambda}$ is nonnegative ($u=1$) and lower semi-continuous with
respect to
the weak topology on $M_1(E)$ ($\hat{L}^{\lambda}u/u$ is bounded). Set
$\tau_n = \inf\{t\geq0 : Y_t\in E_n^c\}$ and
$\tau= \lim_{n\uparrow\infty}\tau_n$. Denote by $\mu_T$ the occupation
measure for $Y$ on $\{T<\tau\}$, which satisfies
\[
\mu_T(B) = \frac{1}{T}\int_0^T1_{Y_t\in B}\,dt\qquad\mbox{for all
$T<\tau
$ and all Borel } B\subseteq E.
\]
For compact $K\subset M_1(E)$, it follows that for all $y\in E$ [\citet
{MR0428471}, Section 7],
%
%
\begin{equation}\label{EDVub}
\limsup_{T\uparrow\infty}\frac{1}{T}\log\hat{P}^{\lambda,y}(\mu
_T\in
K,T<\tau)\leq-\inf_{\mu\in K}I^{\lambda}(\mu).
\end{equation}
Furthermore, the following facts hold.
\begin{lemma}\label{LIfunctprop}
Let Assumptions \ref{regmod}, \ref{regmoddomain} and \ref
{boundedpotential}
hold. If there exists a $\mu^*\in M_1(E)$ such that $I^{\lambda}(\mu^*)
= 0$,
then for all $y\in E$, $\hat{P}^{\lambda,y}(\tau=\infty) = 1$ and
$\mu^*$
possesses a $C^{2,\gamma}(E,\reals)$ density $g^*$, such that
$\tilde{L}^{\lambda}g^* = 0$, where $\tilde{L}^{\lambda}$ is the
formal adjoint
to $\hat{L}^{\lambda}$.
\end{lemma}

\begin{pf}
That
$\hat{P}^{\lambda,y}(\tau= \infty) = 1$ for all $y\in E$ follows by
repeating the argument in {Lemma $3.5$} up to, and through the point
where $\hat{P}^{\lambda,y}(\tau=\infty) = 1$ for $\mu^*$ a.e.
$y\in
E$. By
\citeauthor{MR1326606} [(\citeyear{MR1326606}), Theorem 1.15.1], the same conclusion can be extended
to all
$y\in E$.

As for the second statement, the argument in Lemma $3.6$ can
be repeated
up to and including the point where it is shown that $I^{\lambda}(\mu
^*) =
0$ implies $\int_E \hat{L}^{\lambda}w(x)\mu^*(dx) = 0$ for all
$w\in
C^{\infty}_0(E,\reals)$. Thus the conclusions follow from
\citet{MR1326606}, page 181. In particular, $\int_E g^*\,dy = 1 <
\infty$.
\end{pf}

\begin{pf*}{Proof of Theorem \protect\ref{Tbellmansolution}}
As in the proof of Theorem \ref{TLambdastructure}, the proof of
Theorem \ref{Tbellmansolution} adapts the results in \citet{MR2206349}
to the present setting, now extending {Theorem $3.7$} and
{Theorem $3.8$}.
{Theorem $3.7$} yields that for $\lambda= \lambda_c$, $Y$ is
$(\hat{P}^y)_{y\in E}$-positive recurrent, while {Theorem $3.8$}
states that the solution $v_c$ corresponding to $\lambda_c$ is unique
up to an additive constant. {Theorem $3.8$} carries over with only
notational changes, in the light of \eqref{EAhatArel} and Lemma \ref
{LKSThm26Lemmas}.

{Theorem $3.7$} is composed of three parts. The first part
({Proposition $3.3$}) states that if $Y$ is transient under
$(\hat{P}^{\lambda,y})_{y\in E}$, then there exists $\alpha>0 $ such
that for
all $f\in C_0(E), f\geq0$, $y\in E$, and $T$ large enough
%
%
\begin{equation}\label{EKStranscond}
E_{\hat{P}^{\lambda}}^y[f(Y_T),T<\tau]\leq C(y)e^{-\alpha T}.
\end{equation}
The second part uses the first part to show that $Y$ is recurrent under
$(\hat{P}^{\lambda_c,y})$. The third part states that $Y$ is actually positive
recurrent. The argument used to prove the second part (recurrence) carries
over unchanged, hence, details are provided on the first and third parts.

As for the first part ({Proposition $3.3$}) assume that $Y$ is transient
under $(\hat{P}^{\lambda,y})_{y\in E}$. Consider the set $\mathcal
{C}_m =
\{\mu\in M_1(E) : \mu(E_l)\geq1-\delta_l\ \forall l\geq m\}
$. To
construct the constants $\delta_l$, set
\[
U_0 = -\bigl(V + \tfrac{1}{2}\nabla w'\hat{A}\nabla w + \tilde
{b}'\nabla
w +
\tfrac{1}{2}\trace(AD^2)\bigr)
\]
for the function $w$ from Assumption \ref{generalcase}. Note that by
\eqref{Ewdropoff},
\[
\lim_{n\uparrow\infty}\inf_{y\in E\setminus E_n}
U_0(y) =
\infty.
\]
Let $c$ be as in \eqref{EAhatArel}. Set $\beta_0\,{=}\,\inf_{y\in E}c(U_0(y)\,{+}\,\lambda)$
and \mbox{$\beta_l\,{=}\,\inf_{y\in E\setminus E_l}c(U_0(y)\,{+}\,\lambda)$}. Let $M > 0$ and set
\[
\delta_l = \frac{M + |\beta_0|}{|\beta_0| + \beta_l}.
\]
Since $\beta_l\uparrow\infty$, it follows that $\delta_l\downarrow0$
and hence,
the set $\mathcal{C}_m$ is relatively compact (weak topology) for $m$ large
enough, so that $\beta_m > 0$ and $\delta_m < 1$. Set
$\bar{\phi} = \exp(c(w - v))$ where $v$ is such that $Mv =
\lambda$. Then, repeating the arguments in {Proposition $3.3$},
it follows
that for any $f\in C_0(E)$ [{equation $(3.22)$}],
\[
E_{\hat{P}^{\lambda}}^y[f(Y_T),T<\tau]\leq\sup_{E}|f|\ \hat
{P}^{\lambda,y}[\mu_T\in\mathcal{C}_m,T <\tau] +
\sup_{E}|f/\bar{\phi}|\ \bar{\phi}(y)e^{-M t}.
\]
By \eqref{EDVub} it follows, after taking $M\uparrow\infty$ that
%
%
\begin{equation}\label{EKStranscondlimit}
\limsup_{T\uparrow\infty}\frac{1}{T} \log E_{\hat{P}^{\lambda
}}^y[f(Y_T),T<\tau] \leq
-\inf_{\mu\in\bar{\mathcal{C}}_m} I^{\lambda}(\mu).
\end{equation}
Now,\vspace*{-1pt} since $I^{\lambda}$ is lower semi-continuous and $\bar{\mathcal
{C}}_m$ is
compact, $-\inf_{\mu\in\bar{\mathcal{C}}_m}I^{\lambda}(\mu) =
-I^{\lambda}(\mu^*)$ for some $\mu^*\in\bar{\mathcal{C}}_m$.

It is now shown that $I^{\lambda}(\mu^*) > 0$. Suppose, by
contradiction, that
\mbox{$I^{\lambda}(\mu^*) = 0$} and that $Y$ is transient under $(\hat
{P}^{\lambda,y})_{y\in E}$. Since
$I^{\lambda}(\mu^*) = 0$, Lemma \ref{LIfunctprop} implies that:
\begin{longlist}[(a)]
\item[(a)]
$\hat{P}^{\lambda,y}[\tau=\infty] = 1$ for all $y\in E$.\vadjust{\goodbreak}
\item[(b)]
With $\tilde{L}^{\lambda}$ denoting the adjoint operator to $\hat
{L}^{\lambda}$,
there exists a $C^{2,\gamma}(E,\reals)$ positive function $g^*$ such that
$\tilde{L}^{\lambda}g^* = 0$ and $\int_{E}g^*(x)\,dx = 1$.
\end{longlist}
However, \citeauthor{MR1326606} [(\citeyear{MR1326606}), Corollary 4.9.4] implies that if $Y$ is
transient under
$(\hat{P}^{\lambda,y})_{y\in E}$ and there exists some $\tilde{\phi
}\in
C^{2,\gamma}(E,\reals),\tilde{\phi} > 0$ which satisfies
\mbox{$\tilde{L}^{\lambda}\tilde{\phi} = 0$}, then $\int_E\tilde{\phi
}\,dy <
\infty$
implies $\hat{P}^{\lambda,y}[\tau<\infty] > 0$ for all $y\in E$. This
conclusion contradicts (a) above. Thus, $I^{\lambda}(\mu^*) > 0$ and,
in view of
\eqref{EKStranscondlimit}, the inequality in \eqref{EKStranscond} holds
for any $0 < \alpha< I^{\lambda}(\mu^*)$ and
large enough $T$.

To show the third part [positive recurrence for $Y$ under
$(\hat{P}^{\lambda,y})_{y\in E}$], the same steps as in part one can be
repeated to obtain \eqref{EKStranscondlimit}, as these steps do not
require that $Y$ is transient. Now, if $\inf_{\mu\in\bar{\mathcal
{C}}_m} I^{\lambda}(\mu) > 0$, then for all $f\in C_0(E)$, \eqref
{EKStranscond} implies that
\[
\int_0^\infty E_{\hat{P}^{\lambda}}^y[f(Y_T),T<\tau]\,dT <
\infty,
\]
in which case $Y$ is transient under $(\hat{P}^{\lambda,y})_{y\in E}$
[\citet{MR1326606}, Chapter 4.2]. But part two implies that $Y$ is
recurrent. Thus, there is a $\mu^*\in\bar{\mathcal{C}}_m$ such that
$I^{\lambda}(\mu^*) = 0$. Since $\mu^*$ is a probability measure,
ergodicity follows by Lemma~\ref{LIfunctprop}.
\end{pf*}


\begin{pf*}{Proof of Proposition \protect\ref{Ponestatenice}}
The invariant density for $Y$ under $(\hat{P}^y)_{y\in E}$ is $\tilde
{\phi} = \phi_c^2 m_{\nu}$
where $\phi_c = \exp(v_c/\delta)$ and $m_{\nu}$ is from
\eqref{Emnudef}. Equations \eqref{priest} and~\eqref{duaest} from
Theorem \ref{Tpdemodel} become
\begin{eqnarray*}
E_P^y[(X^{\hatp}_T)^p] &=& e^{\lambda T+v(y)}E_{\hat P}^y[\phi
_c(Y_T)^{-\delta}],\\
E_P^y[(M^{\mfnct}_T)^q]^{1-p} &=& e^{\lambda
T+v(y)}E_{\hat P}^y\bigl[\phi_c(Y_T)^{-{\delta}/{(1-p)}}\bigr]^{1-p}.
\end{eqnarray*}
Since $\int_E m_{\nu}\,dy < \infty$, $\int_E \phi^2_cm_{\nu}\,dy <
\infty$, it
follows that
\[
\int_E\phi_c^{-\delta}\phi_c^2m_{\nu}\,dy < \infty,\qquad
\int_E\phi_c^{-\delta/(1-p)}\phi_c^2m_{\nu}\,dy < \infty
\]
provided that $2-\delta> 0$ and $2-\frac{\delta}{1-p} > 0$. These
conditions are
equivalent to those in \eqref{Eqrhocondnice}. Thus, by the ergodic result
\eqref{Eergodicassumptions} it holds that
\begin{eqnarray*}
\lim_{T\uparrow\infty}E_{\hat P}^y[\phi_c(Y_T)^{-\delta}] &=&
\int_E\phi_c^{-\delta}\phi_c^2m_{\nu}\,dy \equiv K_1,\\
\lim_{T\uparrow\infty}E_{\hat P}^y\bigl[\phi_c(Y_T)^{-{\delta
}/{(1-p)}}\bigr]^{1-p} &=&
\int_E\phi_c^{-\delta/(1-p)}\phi_c^2m_{\nu}\,dy\equiv K_2
\end{eqnarray*}
and long-run optimality follows. Furthermore, in light of \eqref
{celbou}, the
conclusion in~\eqref{ECLEdecaynice} follows with $K =
\frac{1}{p}\log(K_2/K_1)$.
\end{pf*}
%


\section{\texorpdfstring{Proofs of Section \lowercase{\protect\ref{secapp}}}{Proofs of Section 4}}

\begin{pf*}{Proof of Theorem \protect\ref{teolin}}
If $v_0\in\reals^k$ and $v_1\in\reals^{k\times k}$,
$v_1$ symmetric sol\-ve~\eqref{vzeq} and \eqref{vueq},
respectively, then $v(y)=v_0'y-\frac{1}{2}y'v_1 y$ solves
\eqref{EmainvPDE} for $\lambda$ from~\eqref{vleq}. Under Assumption
\ref{assmodlindif}, the condition \eqref{Ewdropoffalt} holds,
hence, Theorems~\ref{TLambdastructure} and~\ref{Tbellmansolution} imply that, if
for $v$
the associated process $Y$ is $\hat{P}^y$-tight in~$\reals^k$ for each
$y\in\reals^k$, then $v$ is the desired solution. The Riccati equation~\eqref{vueq} admits the form
\[
v_1 \mathbf{BB}'v_1 - v_1 \mathbf{A} - \mathbf{A}'v_1 -
\mathbf{C}'\mathbf{C} = 0
\]
with
\[
\mathbf{B} = (A-q\ups'\Sigma^{-1}\ups)^{1/2};\qquad\mathbf{A}
= -(b + q\ups'\Sigma^{-1}\mu_1);\qquad\mathbf{C} =
\sqrt{q}\sigma^{-1}\mu_1,
\]
where $\mathbf{B}$ is assumed to be the unique symmetric positive
definite square root of $A-q\ups'\Sigma^{-1}\ups$. $\mathbf{C}$ is
a real valued matrix when $p<0$. For a real valued square matrix $M$,
write $M>0$ if $M+M'$ is strictly positive definite. If $M > 0$, then
the real part of each of its eigenvalues is strictly positive. To see
this, let $x,\lambda$ such
that $Mx = \lambda x$. Then (with $\bar{x}$ denoting the complex
conjugate of $x$)
\[
0 < \bar{x}'(M+M')x = \bar{x}'\lambda x + \bar{\lambda}\bar{x}'x =
2|x|^2\operatorname{Re}[\lambda].
\]
From \citeauthor{MR1997753} [(\citeyear{MR1997753}), Lemma 2.4.1], if there exist two matrices
$F_1\in
\reals^{k\times
k}$ and $F_2\in\reals^{n\times k}$ such that $\mathbf{A} -
\mathbf{B}F_1 < 0$ and $\mathbf{A}'-\mathbf{C}'F_2 < 0$, then there
exists a unique solution $v_1$ such that $\mathbf{A}-\mathbf{BB}'v_1
< 0$. Since $b>0$ and $p<0$, the choice of $F_1 =
-q\mathbf{B}^{-1}\ups'\Sigma^{-1}\mu_1$ and $F_2 = -\sqrt{q}\rho a'$
suffices. The condition $\mathbf{A}-\mathbf{BB}'v_1 < 0$ yields
%
%
\begin{equation}\label{Enegdefcond}
-\bigl((b + q\ups'\Sigma^{-1}\mu_1)
+(A-q\ups'\Sigma^{-1}\ups)v_1\bigr) < 0.
\end{equation}
Therefore,
\begin{eqnarray*}
&&\bigl(v_1(A-q\ups'\Sigma^{-1}\ups)
+ (b+q\ups'\Sigma^{-1}\mu_1)'\bigr) \\
&&\qquad= \bigl((b +
q\ups'\Sigma^{-1}\mu_1)+
(A-q\ups'\Sigma^{-1}\ups)v_1\bigr)' > 0
\end{eqnarray*}
has eigenvalues with strictly positive real part and is invertible. Thus
$v_0$ from \eqref{vzeq} is well defined.
It remains to prove that $(Y_t)_{t\geq0}$ is $\hat{P}^y$-tight in
$\reals^k$. Under $\hat{P}^y$, $Y$ has the dynamics
\begin{eqnarray*}
dY_t& =& \bigl(-\bigl((b + q\ups'\Sigma^{-1}\mu_1)
+(A-q\ups'\Sigma^{-1}\ups)v_1\bigr)Y_t
\\
&&\hspace*{44pt}{} -q\ups'\Sigma^{-1}\mu_0+
(A-q\ups'\Sigma^{-1}\ups)v_0\bigr)\,dt + a\,dW_t.
\end{eqnarray*}
Setting
\begin{eqnarray*}
\mathbf{D} &=& \bigl((b + q\ups'\Sigma^{-1}\mu_1)
+(A-q\ups'\Sigma^{-1}\ups)v_1\bigr),\\
\mathbf{E}&=&\mathbf{D}^{-1}\bigl(-q\ups'\Sigma^{-1}\mu_0 +
(A-q\ups'\Sigma^{-1}\ups)v_0\bigr),
\end{eqnarray*}
the dynamics takes the form
\[
dY_t = \mathbf{D}(\mathbf{E}-Y_t)\,dt + a\,dW_t,
\]
where $\mathbf{D}\in\reals^{k\times k}$ is such that $\mathbf{D}>0$ and
$\mathbf{E}\in\reals^k$. For $Z_t = Y_t- \mathbf{E}$ it follows that
\[
dZ_t = -\mathbf{D}Z_t \,dt + a\,dW_t.
\]
Since
%
%
\begin{equation}\label{EtempYZcalc}
|Y_t|^2 \leq2|\mathbf{E}|^2 + 2|Z_t|^2,
\end{equation}
is suffices to show that $(Z_t)_{t\geq0}$ is $\hat{P}^y$-tight in
$\reals^k$, which follows because any compact set is contained in a
closed ball around the origin, and by \eqref{EtempYZcalc} if~$Z_t$
is in a closed ball around $0$, then so is $Y_t$. Using the
methods derived in \mbox{\citet{MR0494525}} it follows that under
$(\hat{P}^y)_{y\in E}$, $Z$ is positive recurrent. To show this,
let $\lambda^*$ and $\lambda_{*}$ denote the maximum and minimum
eigenvalues of $A=aa'$ and let $\theta_*$ denote the minimum
eigenvalue of $\mathbf{D}+\mathbf{D}'$. Since $\mathbf{D}>0$ and by
Assumption \ref{assmodlindif}, $A>0$ it follows that $\lambda_*>0$
and $\theta_*>0$. Furthermore,
\[
\inf_{x : |x|=1} x'Ax = \lambda_* \qquad\inf_{x : |x|=1} x'\mathbf{D}x
= \theta_{*}.
\]
Thus, with the notation of \citeauthor{MR1326606} [(\citeyear{MR1326606}), Theorem 6.2],
$\underline{\alpha}(r) = \lambda_{*}$ and
\[
\overline{\beta}(r) = \sup_{x : |x| = 1}\frac{\trace(A) - x'A x - 2r^2
x'\mathbf{D}x}{x'Ax}\leq\frac{\trace(A)-\lambda_{*}-2r^2\theta
_{*}}{\lambda^*},
\]
where the last inequality follows for $r$ large enough so that the
numerator is negative. Therefore, \citeauthor{MR1326606} [(\citeyear{MR1326606}), Theorem 6.2]
applies and $Z$ is positive recurrent.
\end{pf*}


\begin{pf*}{Proof of Proposition \protect\ref{Pgeneralou}}
In light of Theorems \ref{teolin} and \ref{Tlongrunopt}, it suffices
to show that for $p<0$
the following quantity is bounded as a function of $y$:
%
%
\begin{eqnarray}\label{Eboundquant}
&&\bigl(pr_0 - \lambda-
\tfrac{1}{2}q(\nu_0'+b\nu_1'y)(\nu_0 + b\nu_1y)
\nonumber
\\[-8pt]
\\[-8pt]
\nonumber
&&\hspace*{53pt}\qquad{}+ \tfrac{1}{2}q\rho'\rho(v_0-v_1y)^2\bigr)e^{-v_0 y +
({1}/{2})v_1 y^2}.
\end{eqnarray}
Note that, for $p<0$, $\Theta$ in \eqref{EThetaval} satisfies
%
%
\begin{equation}\label{EouThetabound}
\Theta> (1+q\rho'\nu_1)^2.
\end{equation}
Therefore,
\[
v_1 =b\delta\bigl(\sqrt{\Theta} -(1+q\rho'\nu_1)\bigr)>0
\]
and \eqref{Eboundquant} is bounded over $\reals$ only if the
quadratic term is negative
\[
\tfrac{1}{2}q v_1^2\rho'\rho- \tfrac{1}{2}qb^2\nu_1'\nu_1 < 0.
\]
But
\begin{eqnarray*}
\tfrac{1}{2}q v_1^2\rho'\rho- \tfrac{1}{2}qb^2\nu_1'\nu_1 &=&
\tfrac{1}{2}q\rho'\rho\delta^2
b^2\bigl(\sqrt{\Theta}-(1+q\rho'\nu_1)\bigr)^2\\
&&{} -
\tfrac{1}{2}b^2\delta\bigl(\Theta-(1+q\rho'\nu_1)^2\bigr)\\
&=&
\tfrac{1}{2}\delta^2b^2\bigl(\sqrt{\Theta}-(1+q\rho'\nu_1
)\bigr)\\
&&{}\times \bigl(
(2q\rho'\rho-1)\sqrt{\Theta}-(1+q\rho'\nu_1)\bigr).
\end{eqnarray*}
From \eqref{Egeneraloucond}, the quantity
$(2q\rho'\rho-1)\sqrt{\Theta}-(1+q\rho'\nu_1)$
is negative, while\break
$\frac{1}{2}\delta^2b^2(\sqrt{\Theta}-(1+q\rho'\nu_1
))$
is positive by \eqref{EouThetabound}. Therefore, the leading
quadratic term is negative and the result follows.
\end{pf*}


\begin{pf*}{Proof of Corollary \protect\ref{Cfadstypeou}}
When $\nu_1 = -\kappa\rho$ condition \eqref{Egeneraloucond} reduces to
\[
(1-2q\rho'\rho)\bigl(1+q\rho'\rho(\kappa^2-2\kappa
)\bigr)^{1/2}+(1-q\kappa\rho'\rho)
> 0.
\]
Set $x=q\rho'\rho$ and consider the continuous function
\[
f(x,\kappa)
=(1-2x)\bigl(1+x(\kappa^2-2\kappa)
\bigr)^{1/2}+(1-\kappa
x)
\]
on $0<x<1,\kappa\in\reals$. For a fixed $0<x<1$ consider the implicit
equation for~$\kappa$ obtained by setting $f(x,\kappa) = 0$. For
$0<x\leq
\frac{1}{4}$ one can show there are no $\kappa\in\reals$ such that
$f(x,\kappa) = 0$. For $\frac{1}{4}<x<1$, $f(x,\kappa) = 0$ only along
the curve $\kappa= \frac{2}{4x-1}$. For a fixed $x$ and
large positive $\kappa$
\[
f(x,\kappa)\approx\kappa\bigl((1-2x)\sqrt{x}-x\bigr) < 0
\]
and for large negative $\kappa$
\[
f(x,\kappa)\approx|\kappa|\bigl((1-2x)\sqrt{x}+x\bigr) > 0,
\]
therefore,\vspace*{-2pt} plugging back in $q\rho'\rho$ for $x$, for $\frac
{1}{4}<q\rho
'\rho<1$ the restriction $\kappa<
\frac{2}{4q\rho'\rho-1}$ is necessary.
\end{pf*}


\begin{pf*}{Proof of Proposition \protect\ref{PfadsisOK}}
When $\kappa=1$, by Corollary \ref{Cfadstypeou}, long-run
optimality holds for $0 < q\rho'\rho\leq\frac{1}{4}$. For
$\frac{1}{4}<q\rho'\rho< 1$ long-run optimality holds if $1 <
\frac{2}{4q\rho'\rho- 1}$ or $q\rho'\rho< \frac{3}{4}$.
Consider now $q\rho'\rho\geq\frac{3}{4}$, which is equivalent
to $\delta\geq4$ since $\delta= \frac{1}{1-q\rho'\rho}$. When
$\kappa=1$ the solution $v_1,v_0$ and $\Theta$ simplify
considerably to $ \Theta= \delta^{-1}$; $v_1 =
b(\sqrt{\delta}-1)$; and $v_0 = q\delta\rho'\nu_0 $.
Under $\hat{P}$, $Y$ has the dynamics
\[
dY_t = -\frac{b}{\sqrt{\delta}}Y_t \,dt + \,dW_t.
\]

For $Y_0=y$, it follows that $Y_t \sim N(\mu_t,\sigma^2_t)$ with
$\mu_t = ye^{-{b}/{\sqrt{\delta}}t}$, and $\sigma^2_t =
\frac{\sqrt{\delta}}{2b}(1-e^{-2{b}/{\sqrt{\delta}}t})$.
Therefore,
\[
E_{\hat P}^y\bigl[e^{-v(Y_t)}\bigr] = E[\exp(\mathbf{A}Y_t^2 +
\mathbf{B}Y_t)],
\]
where $\mathbf{A}=\frac{b}{2}(\sqrt{\delta}-1)$,
$\mathbf{B}=-q\delta\rho'\nu_0$. For $X\sim N(\mu,\sigma^2)$,
%
%
\begin{eqnarray}\label{Enormalvarcalc}
&& E[e^{\mathbf{A}X^2 + \mathbf{B}X}]
\nonumber\hspace*{-35pt}
\\[-8pt]
\\[-8pt]
\nonumber
&&\qquad=
\cases{
(1-2\mathbf{A}\sigma^2)^{-1/2}\exp\bigl((1-2\mathbf{A}\sigma
^2)^{-1}(\mu^2\mathbf{A} + \mu\mathbf{B} + \frac{1}{2}\sigma^2
\mathbf{B}^2)\bigr),\vspace*{2pt}\cr
 \hspace*{39pt}\mathbf
{A} < \displaystyle\frac{1}{2\sigma^2},\vspace*{5pt}\cr
\infty,\qquad \mathbf{A} \geq\displaystyle\frac{1}{2\sigma^2}.}\hspace*{-35pt}
\end{eqnarray}

Therefore, $E_{\hat P}^y[e^{-v(Y_t)}]<\infty$ if and only if
$\frac{b}{2}(\sqrt{\delta}-1) < \frac{1}{2\sigma^2_t}$.
This condition reduces to
%
%
\begin{equation}\label{Efadsfiniteexprequirement}
1 +
\frac{\sqrt{\delta}}{2}\bigl(1-\sqrt{\delta}\bigr)(1-e^{-
({2b}/{\sqrt{\delta}})t})>
0.
\end{equation}
Note that the left-hand side of \eqref{Efadsfiniteexprequirement} is
equal to $1$ at $t=0$ and monotonically decreasing in $t$ for
$\delta> 1$. Setting this expression equal to $0$, and solving for~$t$ yields
\[
\hat{t} =
-\frac{\sqrt{\delta}}{2b}\log\biggl(\frac{\sqrt{\delta}(\sqrt{\delta
}-1)-2}{\sqrt{\delta}(\sqrt{\delta}-1)}\biggr).
\]
If $\sqrt{\delta}(\sqrt{\delta}-1) > 2$ or equivalently,
$\delta> 4$, then $\hat{t} > 0$ exists. This proves statement~(i)
in Proposition \ref{PfadsisOK}. If $\delta= 4$, the left-hand side of
\eqref{Efadsfiniteexprequirement} reduces to $e^{-b t} > 0$ for
all $t>0$. Thus,
\[
E_{\hat P}^y\bigl[e^{-v(Y_t)}\bigr] = \exp\biggl(\frac{b}{2}(t + y^2) - 4q\rho'\nu_0
e^{{b}/{2} t} y + \frac{8}{b}q^2(\rho'\nu_0)^2(1-e^{-b t})e^{b t}\biggr).
\]
On the other hand, for $\delta= 4$, it follows that
$E_{\hat P}^y[e^{-({1}/{(1-p)})v(Y_t)}]^{1-p}<\infty$ if
\[
\frac{q}{2} + \frac{1}{2(1-p)}e^{-b t} > 0
\]
which is true for all $t>0$ because $\delta= 4$ only when $0 < q<
1$. Thus,
\begin{eqnarray*}
&&E_{\hat P}^y\bigl[e^{-({1}/{(1-p))}v(Y_t)}\bigr]^{1-p} \\
&&\qquad= \bigl(1-(1-q)(1-e^{-b
t})\bigr)^{-{(1-p)}/{2}}\\
&&\qquad\quad{}\times\exp\biggl((1-p)\biggl(\frac{1}{2}b(1-q)y^2e^{-b t} - 4q(1-q)\rho
'\nu_0 e^{-{b}/{2}t}y \\
&&\hspace*{150pt}{}+ \frac{8}{b}q^2(1-q)^2(\rho'\nu
_0)^2(1-e^{-b t})\biggr)\\
&&\hspace*{2pt}\qquad\qquad\qquad{}/{\bigl(1-(1-q)(1-e^{-b t})\bigr)}\biggr).
\end{eqnarray*}
If $\nu_0 = 0$, $E_{\hat P}^y[e^{-v(Y_t)}]\sim Ke^{\frac{b}{2}t}$\vspace*{-1.5pt} and
$E_{\hat P}^y[e^{-({1}/{(1-p)})v(Y_t)}]^{1-p}\sim q^{-{(1-p)}/{2}}$ for
lar\-ge~$t$, so the certainty equivalent loss is bounded by
$-\frac{b}{2p}$, proving (ii) in Proposition~\ref{PfadsisOK}.

If $\nu_0\neq0$, then $E_{\hat P}^y[e^{-v(Y_t)}]\sim e^{K_1e^{b t}}$ and
$E_{\hat P}^y[e^{-({1}/{(1-p)})v(Y_t)}]^{1-p}\sim\allowbreak q^{-{(1-p)}/{2}}e^{K_2}$ for
large $T$, where $K_1, K_2$ are
positive constants. In this case, the certainty equivalent loss
diverges to $-\infty$ with speed of the order of $\frac{K_1}{t}e^{b
t}$. This proves (iii) and the proof is complete.
\end{pf*}


\begin{pf*}{Proof of Proposition \protect\ref{LcirisOK}}
Since $v_0,v_1$ satisfy \eqref{Eciroptans}, $v(y)=v_0\log y + v_1
y$ solves \eqref{EmainvPDE}. Under $\hat{P}$ the dynamics of $Y$
are
\[
dY_t =
\sqrt{\Theta}\biggl(\frac{\sqrt{\Lambda}+({1}/{2})a^2}{\sqrt{\Theta
}}-Y_t\biggr)\,dt
+ a\sqrt{Y_t}\,dW_t
\]
and, as mentioned in Section \ref{SCIR}, the positivity of $\Theta
,\Lambda$
give that $Y$ is $\hat{P}^y$-tight in $(0,\infty)$ for each $y\in
(0,\infty)$.

Therefore, long-run optimality will follow if the quantity $F$ for
$p<0$ from Theorem \ref{Tlongrunopt} is bounded over
$(0,\infty)$. Specifying to this example, it is necessary to show that
%
%
\begin{eqnarray}\label{Ecirbound}
&&\biggl(pr_0 + pr_1 y - \lambda
-\frac{1}{2}q(\nu_0'+y\nu_1)\frac{1}{y}(\nu_0+y\nu
_1)
\nonumber
\\[-9pt]
\\[-9pt]
\nonumber
&&\hspace*{26pt}\qquad{}+ \frac{1}{2}q\biggl(\frac{v_0}{y}+v_1\biggr)a^2\rho'\rho
y\biggl(\frac{v_0}{y} +v_1\biggr)\biggr)e^{-(v_0\log y + v_1
y)}
\end{eqnarray}
is bounded on $y >0$. This expression admits the form
\[
(\mathbf{A} + \mathbf{B}y + \mathbf{C}y^2)y^{-v_0 -1}e^{-v_1 y}.
\]

For $v_0,v_1$ from \eqref{Eciroptans}, by \eqref{Ecompare} it
follows that $v_0 > 0, v_1 < 0$ and so \eqref{Ecirbound} will
follow only if $\mathbf{A} < 0, \mathbf{C} < 0$. As for $\mathbf{A}$,
\begin{eqnarray*}
\mathbf{A} & =& \frac{1}{2}qa^2\rho'\rho v_0^2 - \frac{1}{2}q\nu
_0'\nu
_0\\[-2pt]
&=&
\frac{1}{2}qa^2\rho'\rho\frac{\delta^2}{a^4}\biggl(\sqrt{\Lambda}-
\biggl(b\theta-qa\rho'\nu_0-\frac{1}{2}a^2\biggr)\biggr)^2
\\[-2pt]
&&{}-\frac{1}{2}\frac{\delta}{a^2}\biggl(\Lambda-\biggl(b\theta-qa\rho'\nu
_0-\frac{1}{2}a^2\biggr)^2\biggr)\\[-2pt]
&=&
\frac{1}{2}\frac{\delta^2}{a^2}\biggl(\sqrt{\Lambda}-\biggl(b\theta-qa\rho
'\nu_0-\frac{1}{2}a^2\biggr)\biggr)\\[-2pt]
&&{}\times \biggl((2q\rho'\rho
-1)\sqrt{\Lambda}
-\biggl(b\theta-qa\rho'\nu_0-\frac{1}{2}a^2\biggr)\biggr).
\end{eqnarray*}
From \eqref{parcir}
\[
(2q\rho'\rho-1)\sqrt{\Lambda}
-\bigl(b\theta-qa\rho'\nu_0-\tfrac{1}{2}a^2\bigr) < 0.\vadjust{\goodbreak}
\]
Thus, $\mathbf{A}<0$ since by \eqref{Ecompare}
\[
\frac{1}{2}\frac{\delta^2}{a^2}\biggl(\sqrt{\Lambda}-\biggl(b\theta-qa\rho
'\nu_0-\frac{1}{2}a^2\biggr)\biggr)
> 0.
\]
As for $\mathbf{C}$,
\begin{eqnarray*}
\mathbf{C} &=& \frac{1}{2}qa^2\rho'\rho v_1^2 + pr_1 -
\frac{1}{2}q\nu_1'\nu_1\\
&=&
\frac{1}{2}qa^2\rho'\rho\frac{\delta^2}{a^4}\bigl((b+qa\rho'\nu
_1)-\sqrt{\Theta}\bigr)^2
+ \frac{1}{2}\frac{\delta}{a^2}\bigl((b+qa\rho'\nu_1)^2 -
\Theta\bigr)\\
&=& \frac{1}{2}\frac{\delta^2}{a^2}\bigl((b+qa\rho'\nu_1)
-\sqrt{\Theta}\bigr)\bigl((1-2q\rho'\rho)\sqrt{\Theta} +
(b+qa\rho'\nu_1)\bigr).
\end{eqnarray*}
From \eqref{parcir}
\[
(1-2q\rho'\rho)\sqrt{\Theta} +
(b+qa\rho'\nu_1)
> 0.
\]
Thus, $\mathbf{C}<0$, since by \eqref{Ecompare},
\[
\frac{1}{2}\frac{\delta^2}{a^2}\bigl((b+qa\rho'\nu_1)
-\sqrt{\Theta}\bigr) < 0.
\]
\upqed\end{pf*}

\end{appendix}
\section*{Acknowledgments}
This paper benefited from the helpful comments of seminar
participants at Cornell University, Hitotsubashi University,
University of
Michigan, Princeton University, University of Texas at Austin, the AMS Meeting
in San Diego, the Oberwolfach Workshop on Stochastic Analysis in Finance
and the Sixth Seminar on Stochastic Analysis at Ascona. We are
indebted to an anonymous referee who helped simplify the proof of the
main result and improve its presentation.


%

\printaddresses

\end{document}